\documentclass[11pt,reqno]{amsart}
\usepackage{amsmath,amsfonts,amssymb,amscd,amsthm,amsbsy,epsf}
\usepackage{graphicx,psfrag,comment}

\newtheorem{theorem}{Theorem}

\newtheorem{proposition}[theorem]{Proposition}
\newtheorem{algorithm}[theorem]{Algorithm}
\theoremstyle{remark}
\newtheorem{remark}{Remark}

\theoremstyle{definition}
\newtheorem{defn}{Definition}

\def\bmid{\mathop{\,\big|\,}}
\def\aff{\operatorname{aff}}
\def\Id{\operatorname{Id}}
\def\rank{\operatorname{rank}}
\def\D{{\mathcal D}}
\def\E{{\mathcal E}}

\def\real{{\mathbb R}}
\def\torus{{\mathbb T}}
\def\zed{{\mathbb Z}}

\def\M{{\mathcal M}}

\def\R{{\mathcal R}}

\begin{document}
\title[Efficient algorithms for invariant  tori and their manifolds]
{Computation of   whiskered  invariant tori and their associated manifolds: new fast algorithms}
\author[G. Huguet]{Gemma Huguet}
\address{Centre de Recerca Matem\`atica, Apartat 50, 08193 Bellaterra (Barcelona), Spain}
\email{gemma.huguet@upc.edu}
\author[R. de la Llave]{Rafael de la Llave }
\address{Department of Mathematics, The University of Texas at Austin,
Austin, TX 78712} \email{llave@math.utexas.edu}
\author[Y. Sire]{Yannick Sire}
\address{Universit\'e Paul C\'ezanne, Laboratoire LATP UMR 6632, Marseille, France }
\email{sire@cmi.univ-mrs.fr}
\begin{abstract}
In this paper we present efficient algorithms for the computation  of several invariant
objects for Hamiltonian dynamics. More precisely, we consider KAM tori (i.e
diffeomorphic copies of the torus such that the motion on them is
conjugated to a rigid rotation) both Lagrangian tori (of maximal dimension) and whiskered tori (i.e.
tori with hyperbolic directions which, together with the tangents to
the torus and the symplectic conjugates span the whole tangent
space). In the case of whiskered tori, we also present algorithms to
compute the invariant splitting and  the invariant manifolds
associated to the splitting. We present them both for the case of discrete time and for differential equations.

The algorithms for tori are based on a Newton method to solve an
appropriately chosen functional equation that expresses invariance.
Among their features we highlight:
\begin{itemize}
\item
The algorithms are efficient: if we discretize the objects by $N$
elements, one step of the Newton method requires only $O(N)$ storage
and $O(N \ln(N))$ operations. Furthermore, if the object we consider
is of dimension $\ell$, we only need to compute functions of $\ell$
variables, independently of what is the dimension of the phase
space.
\item
The algorithms do not require that the system is presented in
action-angle variables nor that it is close to integrable.
\item
The algorithms are backed up by rigorous \emph{a-posteriori} bounds which
state that if the equations are solved with a small residual and some
explicitly computable
condition numbers are not too big, then, there is a true solution which is
close to the computed one.
\item
The algorithms  apply both to primary (i.e non-contractible) and
secondary tori (i.e. contractible to a torus of lower 
dimension, such as islands). They also apply to whiskered tori. 
\end{itemize}
The algorithms for invariant splittings are based on 
computations of proyections (rather than in graph transforms). 
The computations of invariant manifolds are also efficient in 
the sense indicated before.

The algorithms we present have already been implemented. We will report on the
technicalities of the implementation and the results of running them
elsewhere.

\end{abstract}


\maketitle

\section{Introduction}
\label{sec:intro} The goal of this paper is to present efficient
algorithms to compute very accurately several objects of interest in
Hamiltonian dynamical systems (both discrete-time dynamical systems
and differential equations). More precisely, we present algorithms
to compute:
\begin{itemize}
\item
Lagrangian KAM  tori.
\item
Whiskered KAM tori.
\item
The invariant bundles of the whiskered tori.
\item
The stable and unstable manifolds of the whiskered tori.
\end{itemize}
The algorithms are very different. For example, the algorithms for 
tori require the use of small divisors and symplectic geometry
and the algorithms for invariant bundles and invariant manifolds 
rely on the theory of normal hyperbolicity and dichotomies. The computation of 
whikered tori has to combine both.

We recall that KAM tori are manifolds diffeomorphic to a torus which
are invariant for a map or flow, on which the motion of the system
is conjugate to a rotation. As we will see later, this is also
equivalent to quasi-periodic solutions. The tori are called
Lagrangian when they are Lagrangian manifolds, which in our case
amounts just to the fact that the tori have a dimension equal to the
number of degrees of freedom of the system. The tori are called
whiskered when the linearized equation has directions that decrease
exponentially either in the future (stable) or in the past
(unstable)  and these directions together with the tangent to the
torus and its symplectic conjugate span the whole tangent space.
These invariant spaces for the linearization have non-linear
analogues, namely invariant manifolds. It has been recognized since
\cite{Arnold64} that whiskered tori and their invariant manifolds
are very interesting landmarks that organize the long-term behavior
of many systems.

The algorithms we present are based on running an \emph{efficient}
 Newton method to
solve a functional equation, which expresses the dynamical
properties above. What we mean by \emph{efficient} is that if we
discretize the problem using $N$ Fourier coefficients, we require
$O(N)$ storage and only $O(N \ln(N))$ operations for the Newton
step. Since the functions we are  considering are analytic, we see
that the truncation error is $O(\exp(-C N^{1/d} ))$ where $d$ is the
dimension of the object. Note that, in contrast, a straightforward
implementation of a Newton method would require to use $O(N^2)$
storage -- to store the linearization matrix and its inverse -- and
$O(N^3)$ operations to invert.

In practical applications, using the algorithms described in this
paper, computing with several million coefficients becomes quite
practical in a typical desktop computer of today. Implementation
details and the results of several runs will be discussed in another
companion paper \cite{HuguetLS10c}. Given the characteristics of
today's computers, savings in storage space are more crucial than
savings in operations for these problems.

The algorithms we present here are inspired by the rigorous results
of \cite{LlaveGJV05} -- for KAM tori -- and
\cite{FontichLS09a,FontichLS09b} --for whiskered tori. The algorithms to compute the stable and unstable manifolds had not been previously 
discussed. The rigorous
results of the above papers are also based on a  Newton method
applied to the same functional equation that we consider here.
Of course, going from a mathematical treatment to 
a practical algorithm requires significantly many more details. 
In particular, the algorithms to compute invariant 
splittings and invariant manifolds are different 
from those in the above references. 
This paper discusses these algorithmic issues. 

The results of the papers \cite{LlaveGJV05,FontichLS09b} give a
justification of the algorithms for tori and splittings 
presented here. The theorems in
\cite{LlaveGJV05,FontichLS09b}, have been formulated in an
\emph{a-posteriori} way, i.e. the theorems assert
that  if we have a function which solves
approximately the invariance equation very accurately (e.g. the
outcome of a successful run of the algorithms) and which also
satisfies some explicit non-degeneracy conditions, then, we can
conclude that there is a true solution which is close to the
computed solution. Hence, by supplementing our calculations with the
(very simple) computations of the non-degeneracy conditions (they
play a role very similar to the condition numbers common in
numerical analysis), we can be sure that the computation that we are
performing is meaningful. This allows to compute with confidence
even close to the limit of validity of the KAM theorem (a rather
delicate boundary since the smooth KAM tori do not disappear
completely but rather morph into Cantor sets).

Since the papers \cite{LlaveGJV05,FontichLS09b} contain estimates,
in the present paper, we will only discuss the algorithmic issues.
For example, we will detail how solutions of equations (whose
existence was shown in the above papers) can be computed with small
requirements of storage and small operation count. Note that
different algorithms of a same mathematical operation can have
widely different operation counts and storage requirements. (See,
for example, the discussion in \cite{Knuth97} on the different
algorithms to multiply matrices, polynomials, etc.) On the other
hand, we will not include some implementation issues (methods of
storage of arrays, ordering of loops, precision, etc.) needed to
obtain actual results in a real computer. They will be given in
another paper together with experimental results obtained by 
running the algorithms.

One remarkable feature of the algorithms presented  here is that
they do not require the system to be close to integrable. We only
need a good initial guess for the Newton method. Typically, one uses
a continuation method starting from an integrable case, where
solutions can be computed analytically. However, in the case of
secondary KAM tori, which do not exist in the integrable case, 
one can use, for instance, Lindstedt series, variational methods or 
approximation by periodic orbits to obtain an initial guess.

As for the algorithms to compute invariant splittings, we 
depart from the standard mathematical methods (most of the time based
on graph transforms) and we have found more efficient to device 
an equation for the invariant projections. We also found some 
acceleration of convergence methods that give superexponential  convergence. 
They are based on fast algorithms to solve cohomology equations 
which could be of independent interest (See Appendix~\ref{sec:cohomology}). 

The algorithms to compute invariant manifolds are based on 
the parameterization method \cite{CabreFL03a, CabreFL05}. 
Compared to standard methods such as the graph transform it has 
the advantage that to compute geometric objects of 
dimension $\ell$, we only need to compute with functions of 
dimension $\ell$. In contrast with  \cite{CabreFL03a, CabreFL05}, 
which was based on contractive iterations,
our method is based on a Newton iteration which we also 
implement without requiring large matrices and requiring only 
$N\log(N)$ operations.

\subsection*{An overview of the method}
\label{sec:overview}

The numerical method we use is based on the parameterization methods
introduced in \cite{CabreFL03a,CabreFL03b}. In this
section, we provide a sketch of the issues, postponing some
important details.  For brevity, we make the presentation for
maps only, even if a similar sketch can be made for flows. 

\subsubsection*{Invariant tori}

We observe that if $F$ is a map and we can find an embedding
$K$ in which the motion on the torus is a rotation $\omega$, it
should satisfy  the equation
\begin{equation}
\label{eq:inv1}
F( K( \theta) ) - K(\theta + \omega)  = 0.
\end{equation}

Given an approximate solution of
\eqref{eq:inv1}, i.e.
$$F( K( \theta) ) -  K(\theta + \omega)  = E(\theta),$$
the Newton method aims to find $\Delta$ solving
\begin{equation}
\label{eq:correction} D F( K( \theta) )\Delta(\theta) -
\Delta(\theta + \omega)  = - E(\theta),
\end{equation}
so that $K+\Delta$ will be a much more approximate solution.

The main idea of the Newton method is that, using the decomposition
into invariant subspaces, one can decompose
\eqref{eq:correction} into three components
\begin{equation}
\label{eq:correction_split}
\begin{split}
& D F( K( \theta) )\Delta^s(\theta) - \Delta^s(\theta + \omega)
= -E^s(\theta)\\
& D F( K( \theta) )\Delta^u (\theta) - \Delta^u(\theta + \omega)
= -E^u(\theta)\\
& D F( K( \theta) )\Delta^c (\theta) - \Delta^c(\theta + \omega)
= - E^c (\theta)\\
\end{split}
\end{equation}
where the $s,u$ refer to the stable, unstable components
and $c$ refers to the component along the tangent to the torus
and its symplectic conjugate. For Lagrangian tori, only the $E^c$ part appears in the equations.

The algorithm requires:
\begin{itemize}
\item
Efficient  methods to evaluate the LHS of \eqref{eq:inv1}.
\item
Efficient  methods to compute the splitting.
\item
Efficient methods to solve the equations
\eqref{eq:correction_split}.
\end{itemize}

As we will see in Section \ref{sec:FT}, to evaluate \eqref{eq:inv1},
it is efficient to use both a Fourier representation (which makes
easy to evaluate $K(\theta + \omega)$) and a real representation
which makes easy to evaluate $F( K(\theta))$. Of course, both of
them are linked through the Fast Fourier Transform (FFT from now on).

The methods to compute the splitting are discussed in
Section~\ref{sec:projections}. More precisely, we present a
numerical procedure to compute the projections on the linear
stable/unstables subspaces based on a Newton method. In
\cite{HuguetLS10c}, we present an alternative procedure for the
computation of the projections based on the calculation of invariant
bundles for cocycles. Indeed, these algorithms require the
computation of the projections over the linear subspaces of the
linear cocycle.

The solution of the hyperbolic components in equation
\eqref{eq:correction_split} is discussed in Section~\ref{generalst}
and Appendix~\ref{sec:cohomology}. Indeed, equations of this form
appear as well in the calculation of the invariant splitting
discussed in Section~\ref{sec:projections}. A first method is based
on an acceleration of the fixed point iteration
(Appendix~\ref{sec:cohomiter}). We note that to obtain
superexponential convergence for the solution of
\eqref{eq:correction_split}, we need to use both the Fourier
representation and the real space representation.

In the case that the bundles are one-dimensional, there is yet
another algorithm, which is even faster than the previous ones (see
Appendix~\ref{sec:cohom1d}). The algorithms are discussed for maps,
and they do not have an easy analog for flows except by passing for 
the integration of differential equations.
We think that this is one case where working with time-$1$ maps
is advantageous. 

The most challenging step is the solution of  the center component
of \eqref{eq:correction_split}. This depends on cancelations which
use the symplectic structure, involves small divisors and requires
that certain 
obstructions vanish. Using several geometric identities that take advantage
of the fact that the map is symplectic  (see Section
\ref{fast-whisktori}), the solution of \eqref{eq:correction_split}
in the center direction is reduced to solving the following equation  for
$\varphi$ given $\eta$,
\begin{equation} \label{eq:deltaerr}
\varphi(\theta) - \varphi(\theta + \omega) = \eta(\theta).
\end{equation}

Equation \eqref{eq:deltaerr} can be readily solved using Fourier
coefficients provided that $\int \eta =  0$ (and that $\omega$ is
sufficiently irrational). The solution is unique up to addition of a
constant.

The existence of obstructions -- which are finite dimensional -- is
one of the main complications of the problem. It is possible to show
that, when the map $F$ is exact symplectic, the obstructions for the
solution are $O(||E||^2)$. An alternative method to deal with these
obstructions is to add some new -- finite dimensional -- unknowns
$\lambda$  and solve, instead of \eqref{eq:inv1}, the equation
\[
F(K(\theta) - K(\theta + \omega) + G(\theta)\lambda = 0
\]
where $G(\theta)$ is an explicit function. Even if $\lambda$ is kept
through the iteration involving approximate solutions, it can be
shown that, if the map is exact symplectic, we have $\lambda = 0$.
This counterterm approach also helps to weaken non-degeneracy
assumptions.

A minor issue that we omit in this preliminary discussion is that
the solutions of \eqref{eq:inv1} are not unique. If $K$ is a
solution, $\tilde K$ defined by $\tilde K(\theta) = K(\theta
+\sigma)$ is also a solution for any $\sigma \in \mathbb{R}^{\ell}$.
This can be easily solved by taking an appropriate normalization
that fixes the origin of coordinates in the torus. In
\cite{FontichLS09b} it is shown that this is the only non-uniqueness
phenomenon of the equation. Furthermore, this local uniqueness
property allows to deduce results for vector fields from the results
for maps.

It is important to remark that the algorithms that we will present
can compute  in a unified way both primary and secondary tori. We
recall here that \emph{secondary tori} are invariant tori which are
contractible to a torus of lower dimension, whereas this is not the
case for primary tori. The tori which appear in integrable systems
in action-angle variables are always primary. In quasi-integrable
systems, the tori which appear through Lindstedt series or other
perturbative expansions starting from those of the integrable system
are always primary. Secondary tori, however, are generated by
resonances. In numerical explorations, secondary tori are very
prominent features that have been called  ``islands''. In
\cite{HaroL00}, one can find arguments showing that these solutions
are very abundant in systems of coupled oscillators. As an example
of the importance of secondary tori we will mention that in the
recent paper \cite{DelshamsLS06} they constituted the essential
object to overcome the ``large gap problem'' and prove the existence
of diffusion. In \cite{DelshamsH09}, one can find a detailed
analysis of secondary tori.

In this paper, we will mainly discuss algorithms for systems with
dynamics described by diffeomorphisms. For systems described through
vector fields, we note that, taking time$-1$ maps or, more
efficiently, surfaces of section, we can reduce the problem with
vector fields to a problem with diffeomorphisms. However, in some
practical applications, it may be convenient to have a direct
treatment of the system described by vector fields. For this reason,
we have included the invariance equations for flows, in parallel
with the invariance equations for maps and we have left for the
Appendix the algorithms that are specially designed for flows.

\subsubsection*{Invariant manifolds attached to invariant tori}

When the torus is whiskered, it has invariant manifolds attached to
it. For simplicity, in this presentation we will discuss the case of
one dimensional directions -- even if the torus can be of higher
dimension.

We use again a parameterization method. Consider an embedding $W$ in
which the motion on the torus is a rotation $\omega$ and the motion
on the stable (unstable) whisker consists of a contraction
(expansion) at rate $\mu$, it should satisfy the invariance equation
\begin{equation}
\label{eq:whiskers1} F( W(\theta, s) ) - W(\theta + \omega, \mu s)
= 0.
\end{equation}

Again, the key point is that taking advantage of the geometry of the
problem we can devise algorithms which implement a Newton step to
solve equation \eqref{eq:whiskers1} without having to store---and
much less invert---a large matrix.  We first discuss the so-called
order by order method, which serves as a comparison with more
efficient methods based on the reducibility. Although they are based on the same idea as for the case of tori, they have not been introduced previously and 
constitute one of the main novelties of this paper. We present algorithms that compute
at the same time the torus and the whiskers and algorithms that
given a torus and the linear space compute the invariant manifold
tangent to it. It is clearly possible to extend the method to
compute stable and unstable manifolds in general dimensions (or even
non-resonant bundles). To avoid increasing  the length of this paper
and since higher dimensional examples are harder numerically, we
postpone this to a future paper.


\subsection*{Some remarks on the literature}
\label{sec:literature}
Invariant tori in Hamiltonian dynamics have been recognized as important
landmarks in Hamiltonian dynamics. In the case of whiskered tori,
their manifolds have also been crucial for the study of Arnold diffusion.

Since the mathematical literature is so vast, we cannot hope
to summarize it here. We refer to the rather extensive references
of \cite{Llave01} for Lagrangian tori and those of
\cite{FontichLS09b} for whiskered tori. We
will just briefly mention that \cite{Graff74, Zehnder76}
the earliest references on 
whiskered theory, 
as well as most  of the later references, 
are based on \emph{transformation theory}, that is
making changes of variables that reduce the perturbed Hamiltonian
to a simple form which obviously presents the invariant torus.
{From} the point of view of numerics, this has the disadvantage
that transformations are very hard to implement.

The numerical literature is not as broad as the rigorous one, but it
is still quite extensive. The papers closest to our problems are
\cite{HaroL06a,HaroL06b,HaroL07}, which consider tori of systems
under quasi-periodic perturbation.  
 These papers also contain a
rather wide bibliography on papers devoted to numerical computation
of invariant circles.  Among the papers not included in the
references of the papers above because these appeared later, we
mention \cite{CallejaL08}, which presents other algorithms which
apply to variational problems (even if they do not have a
Hamiltonian interpretation). Another fast method  is the
\emph{``fractional iteration method''} \cite{Simo00}.
Note that the problems considerd in 
\cite{HaroL06a,HaroL06b,HaroL07} do not involve center directions 
(and hence, do not deal with small divisors) and that the frequency 
and one of the coordinates of the torus is given by the external 
perturbation. The methods of \cite{HaroL06a,HaroL06b,HaroL07}
work even if the system is not symplectic (even if they can 
take advantage of the symplectic structure).

The papers  \cite{JorbaO03, JorbaO09} present and implement
calculations of \emph{reducible} tori. This includes tori with
normally elliptic directions.  The use of reducibility indeed leads
to very fast Newton steps, but it still requires the storage of a
large matrix for the changes of variables. As seen in the examples
in \cite{HaroL07,HaroL06c}, reducibility may fail in a codimension
$1$ set (a Cantor set of codimension $1$ manifolds for elliptic 
tori in Hamiltonian systems). For these
reasons, we will not discuss methods based on reducibility in this
paper (even if it is a useful and practical 
tool) and just refer to the references just mentioned. Indeed,
thanks to hyperbolicity, reducibility is not needed in the present
paper.

The paper is organized as follows. In Section~\ref{sec:conventions}
we summarize the notions of mechanics and symplectic geometry we
will use. In Section~\ref{sec:equations} we formulate the invariance
equations for the objects of interest (invariant tori, invariant
bundles and invariant manifolds) and we will present some
generalities about the numerical algorithms.

Algorithms for whiskered tori are discussed in
Section~\ref{sec:whiskered}. In particular, we discuss how to
compute the decomposition \eqref{eq:correction_split} of the
linearized equation \eqref{eq:correction}, and how to solve
efficiently each equation in \eqref{eq:correction_split}.

In Section~\ref{sec:rank1im} we discuss fast algorithms to compute
rank-1 (un)stable manifolds of whiskered tori. More precisely, we
present an efficient Newton method to solve equation
\eqref{eq:whiskers1}.

In Appendix~\ref{sec:cohomology} one can find the fast algorithms to
solve cohomology equations with non-constant coefficients that will
be used in the computation of the splitting
\eqref{eq:correction_split} as well as to solve the hyperbolic
components of equations \eqref{eq:correction_split}. In Appendices
\ref{ap:a}-\ref{ap:d}, one can find the algorithms specially
designed for flows, analogous to the ones for maps.

\section{Setup and conventions}\label{sec:conventions}

We will be working with systems defined on an Euclidean phase space
endowed with a symplectic structure. The phase space under consideration will be

$$\mathcal M \subset \real^{2d-\ell} \times \torus^\ell.$$

We do not assume that the coordinates in the phase space are
action-angle variables. Indeed, there
are several systems (even quasi-integrable ones) which are very
smooth in Cartesian coordinates but less smooth in action-angle
variables (e.g., neighborhoods of elliptic fixed points
\cite{FassoGB98a,GuzzoFB98b}, hydrogen
atoms in crossed electric and magnetic fields
 \cite{RakovicChu95,RakovicChu97} and
several problems in celestial mechanics \cite{CellettiC07}).

We will assume that the Euclidean manifold $\mathcal M$ is endowed with an exact symplectic structure $\Omega=d\alpha$  (for some one-form $\alpha $) and we have
\begin{equation*}
\Omega_z (u,v) = \langle u,J(z)v\rangle,
\end{equation*}
where $\langle \cdot, \cdot \rangle$ denotes the inner product on the tangent space of $\mathcal M$ and $J(z)$ is a skew-symmetric matrix.

An important particular case is when $J$ induces an almost-complex structure, i.e.
\begin{equation}\label{complex}
J^2 = -\Id.
\end{equation}
Most of our calculations do not need this assumption. One
important case, where the identity  \eqref{complex} is not
satisfied, is when $J$ is a symplectic structure on surfaces of section chosen arbitrarily in the
energy surface or when $J$ is the symplectic form expressed in
symplectic polar coordinates near an elliptic fixed 
point.  When \eqref{complex} holds, some calculations can be
made faster.

As previously mentioned, we will be considering
systems described either  by diffeomorphisms or by vector-fields.

\subsection{Systems described by diffeomorphisms}
We will consider maps $F: \mathcal U \subset \mathcal M \mapsto \mathcal M$
which are not only symplectic (i.e. $F^{*} \Omega=\Omega$) but exact symplectic,
 that is
\begin{equation*}
F^* \alpha = \alpha +d P,
\end{equation*}
for some smooth function $P$, called the \emph{primitive function}.

We will also need Diophantine properties for the frequencies of
the torus. For the case of maps, the  useful notion of a Diophantine frequency is:
\begin{equation*}
 \D(\nu,\tau) = \left\{ \omega\in\real^\ell \bmid
|\omega \cdot k-n|^{-1} \le \nu |k|^\tau\ \forall\ k\in \zed^\ell -
\{0\},\ n\in\zed\right\},\,\,\,\nu >\ell.
\end{equation*}

\subsection{Systems described by vector fields}
We will assume that the system is described by a globally
Hamiltonian vector-field $X$, that is
\begin{equation*}
X = J\nabla H,
\end{equation*}
where $H$ is a globally defined function on $T^*\mathcal M$.

In the case of flows, the appropriate notion of
Diophantine numbers is:
\begin{equation*}
\D^{\aff} (\nu,\tau) =  \left\{ \omega\in \real^\ell \bmid |\omega
\cdot k|^{-1} \leq \nu |k|^\tau
\ \forall\ k\in \zed^\ell - \{0\} \right\},\,\,\,\nu \geq \ell -1\\
\end{equation*}

\begin{remark}
It is well known that for non-Diophantine frequencies substantially
complicated behavior can appear \cite{Herman92,FayadKW01}. Observing
convincingly Liouvillian behaviors seems a very ambitious
challenge for numerical exploration.
\end{remark}

\section{Equations for invariance} \label{sec:equations}

In this section, we discuss the functional equations for the
objects of interest, that is, the invariant tori and the associated
whiskers. These functional equations, which describe the invariance of the objects under consideration, are
the cornerstone of the algorithms.

\subsection{Functional equations for whiskered invariant tori for diffeomorphisms}

At least at the formal level, it is natural to search quasi-periodic solutions with frequency $\omega$ (independent over the integers) under the form of Fourier series

\begin{equation}
\label{representation}
x^{(n)} = \sum_{k\in\zed^\ell} \hat x_k e^{2\pi
ik\cdot \omega n}\ ,
\end{equation}
where $\omega \in \real^{\ell}$ and $n \in \zed$.

We allow some components of $x$ in 
\eqref{representation} to be angles. In that
case, it suffices to take some of the components of $x$ modulo~1.

It is then natural to describe a quasi-periodic function using the so-called
``hull'' function $K :\torus^\ell \to \mathcal M$
defined by
\begin{equation*}
K(\theta) = \sum_{k\in\zed^{\ell}} \hat x_k e^{2\pi ik \cdot \theta},
\end{equation*}
so that we can write
\begin{equation*}
x^{(n)} = K(n\omega).
\end{equation*}

The geometric interpretation of the hull function is that it gives
an embedding from $\torus^\ell$ into the phase space. In our
applications, the embedding will actually be an immersion.

It is clear that quasi-periodic functions will be orbits for a map $F$ if and only if the hull function $K$
satisfies:
\begin{equation}\label{invariancemaps}
F\circ K - K\circ T_\omega = 0,
\end{equation}
where $T_{\omega}$ denotes a rigid rotation
\begin{equation}\label{omegatranslation}
T_\omega (\theta) = \theta + \omega.
\end{equation}

A modification of the invariance equations \eqref{invariancemaps} which
we will be important for our purpose consists in considering
\begin{equation}\label{invariance-modifiedmaps}
F\circ K - K\circ T_\omega - (J (K_0)^{-1} DK_0)
\circ T_\omega \lambda =0,
\end{equation}
where the unknowns are now $K:\torus^\ell \to \mathcal M$ (as before) and $\lambda \in \real^\ell$.
Here, $K_0$ denotes a given approximate (in a suitable sense which will be given below) solution of the equation \eqref{invariancemaps}.

It has been shown in \cite{FontichLS09a, FontichLS09b}
(the \emph{vanishing lemma}  that, for exact symplectic
maps, if $(K,\lambda)$ satisfy the equation
\eqref{invariance-modifiedmaps} with $K_0$ close to $K$, then at the end of the iteration of the Newton method, we have
$\lambda=0$ and, therefore,  $K$ is a solution of the invariance equation
\eqref{invariancemaps}.
In other words, the formulations \eqref{invariance-modifiedmaps}
and \eqref{invariancemaps} 
are equivalent.
Of course, for approximate solutions of the
invariance equation \eqref{invariancemaps}, there is no reason why
$\lambda$ should vanish
and it is numerically advantageous to solve 
the equation with more variables.   


The advantage of equation \eqref{invariance-modifiedmaps} for
numerical calculations is that, at the initial stages of the method,
when the error in the invariance equation is large, it is not easy
to ensure that certain compatibility conditions
(see \eqref{eq:compatibility} in Section \ref{cohomology}), are
satisfied approximately, so that the standard Newton method has
problems proceeding. On the other hand, we can always proceed by
adjusting the $\lambda$. This is particularly important for the case
of secondary tori that we will discuss in
Section~\ref{sec:secondary}. We also note that this procedure makes 
possible to deal with tori when the twist condition degenerates.

The equations \eqref{invariancemaps} and \eqref{invariance-modifiedmaps}
will be the centerpiece of our treatment. We will discretize them
using Fourier series and study numerical methods to solve the
discretized equations.

It is important to remark that there are \emph{a posteriori}  theorems
(see  \cite{FontichLS09a, FontichLS09b})
for equations \eqref{invariancemaps}, \eqref{invariance-modifiedmaps}
(as well as their analogous for flows \eqref{invariancevf},
\eqref{invariance-modifiedvf} ). That is,
theorems that ensure that given a function that satisfies
\eqref{invariancemaps}, \eqref{invariance-modifiedmaps}
up to a small error and that, at the same time,
satisfies some non-degeneracy conditions (which are
given quite explicitly),  then there is a true
solution close to the computed one.
Hence, if we monitor the non-degeneracy conditions, we can be sure
that the computed solutions correspond to some real effects and are
not spurious solutions.

\begin{remark}
Notice that for whiskered tori the dimension of the torus $\ell$ is
smaller than half the dimension of the phase space $2d$. Hence, the
algorithms presented  here have the advantage that  they look for a
function $K$  of $\ell$ variables 
to compute invariant objects of dimension $\ell$. This is
important because the cost of handling functions grows exponentially
fast with the number of variables. Indeed, to discretize
a function of $\ell$ variables in a grid of side $h$ into $\real^{2d}$, one needs
to store $(1/h)^\ell\cdot 2d$ real values.
\end{remark}

\begin{remark}\label{nonunique}
Equations \eqref{invariancemaps} and \eqref{invariance-modifiedmaps}
do not have unique solutions. Observe that if $K$ is a solution, for
any $\sigma \in \real^\ell$, $K\circ T_\sigma$ is also a solution.
In \cite{FontichLS09b}, it is shown that, in many circumstances,
this is the only non uniqueness phenomenon in a sufficiently small
neighborhood of $K$. Hence, it is easy to get rid of it by imposing
some normalization. See Section~\ref{sec:uniqueness} for a discussion on this issue.
\end{remark}

\subsection{Functional equations for whiskered invariant tori for vector-fields}

In this case, one can write
\begin{equation*}
x(t) = \sum_{k\in\zed^\ell} \hat x_k e^{2\pi ik\cdot \omega t}
\end{equation*}
where $\omega \in \real^{\ell}$, $t \in \real$ and then the hull function $K$ is defined by
$$x(t)=K(\omega t). $$

The invariance equation for flows is written:
\begin{equation}\label{invariancevf}
\partial_\omega K - X \circ K= 0,
\end{equation}
where $\partial_\omega$ denotes the derivative in direction $\omega$
\begin{equation}\label{eq:deromega}
\partial_\omega= \sum_{k=1}^\ell \omega_k \partial_{\theta_k}.
\end{equation}

The modification of \eqref{invariancevf} incorporating a counterterm
is:
\begin{equation}\label{invariance-modifiedvf}
\partial_\omega K - X\circ K - J (K_0)^{-1} (DX\circ K_0) \lambda
=0,
\end{equation}
where $K_0$ is a given embedding satisfying some non-degeneracy
conditions.

\begin{remark}\label{timeone}
Recall that, taking time$-1$ maps, one can reduce the problem of
vector fields to the problem of diffeomorphisms. Furthermore, since
autonomous Hamiltonian systems preserve energy, we can take a
surface of section and deal with the return map. This reduces by $2$
the dimension of the phase space and the parameterization of the
torus requires $1$ variable less. In practice, it is much more
efficient to use a numerical integrator to compute the point of
intersection with the surface of section than to deal with functions
of one more variable and with two more components.
\end{remark}

\subsection{Some global topological considerations}
\label{sec:global}

In our context, both the domain $\torus^\ell$ and the range of $K$
have topology. As a consequence, there will be some topological
considerations in the way that the torus $\torus^\ell$ gets embedded
in the phase space. More explicitly, the angle variables of
$\torus^\ell$ can get wrapped around in different ways in the phase
space.

A concise way of characterizing the topology of the embedding is to
consider the lift of $K$ to the universal cover, i.e.
\begin{equation*}
\widehat K :\real^\ell \to \real^{2d-\ell} \times\real^\ell,
\end{equation*}
in such a way that $K$ is obtained from $\widehat K$ by identifying
variables in the domain and in the range that differ by an integer.

It is therefore clear that $\forall\ e\in \zed^\ell$
\begin{equation}\label{lifting}
\begin{split}
\widehat K_p (\theta +e) & = \widehat K_p (\theta),\\
\noalign{\vskip6pt} \widehat K_q (\theta +e) & = \widehat K_q
(\theta) + I(e),
\end{split}
\end{equation}
where $\widehat K_p,\widehat K_q $ denote the projections of the
lift on the $p$ and $q$ coordinates of $\real^{2d-\ell}
\times\real^\ell$. It is easy to see that $I(e)$ is a linear
function of $e$, namely
\begin{equation}\label{eq:I}
I(e)_{i=1,\ldots,\ell} = \bigg( \sum_{j=1}^\ell I_{ij}
e_j\bigg)_{i=1,\ldots,\ell}
\end{equation}
with $I_{ij}\in\zed$.

We note that if a function $\widehat K_q$ satisfies
\begin{equation*}
\widehat K_q (\theta +e) = \widehat K_q(\theta) + I(e)\ ,
\end{equation*}
the function
\begin{equation}\label{eq:Kperiodic}
\widetilde K_q(\theta) \equiv \widehat K_q(\theta) - I(\theta)
\end{equation}
is $e-$periodic. The numerical methods will always be based on studying
the periodic functions $\widetilde K_q$, but we will not emphasize
this unless it can lead to confusion.

Of course, the integer valued matrix $I=\left \{I_{ij} \right \}_{ij}$ remains constant if we
modify the embedding slightly. Hence, it remains constant under
continuous deformation. For example, in the integrable case with
$\ell=d$, invariant tori satisfy $\widehat K_q (\theta)=\theta$, so
that we have $I = \Id$. Hence, all the invariant tori
which can be continued from tori of the integrable system will also have
$I=\Id$.

\subsection{Secondary tori}\label{sec:secondary}

One can produce other $\ell$-dimensional tori for which
the range of $I$ is of dimension less than $\ell$.
These tori are known as
\emph{secondary tori}.
 It is easy to see that if
$\rank(I)<\ell$ we can contract $K(\torus^\ell)$ to a diffeomorphic
copy of $\torus^{\rank(I)}$. Even in the case  of maximal tori
$\ell=d$, one can have contractible directions. The most famous
example of this phenomenon are the ``islands'' generated in twist
maps around resonances.

Secondary tori do  not exist in the integrable system and they
cannot be even continuously deformed into some of the tori presented
in the integrable system. This is often described informally as
saying that the secondary tori are generated by the resonances.

Perturbative proofs of existence of secondary tori are done in
\cite{LlaveW04} and in \cite{DelshamsLS06} and in more detail in
\cite{DelshamsH09}. In \cite{Duarte94} one can find rigorous results
showing that these islands have to be rather abundant (in different
precise meanings) in many classes of 2D-maps. In particular, for
standard-like maps, secondary tori appear at arbitrarily large values
of the parameter.

In \cite{HaroL00}, there are heuristic arguments and numerical
simulations arguing that in systems of coupled oscillators, the tori
with contractible directions are much more abundant than the
invariant tori which can be continued from the integrable limit.

In view of these reasons, we will pay special attention to the
computation of secondary tori.

We want to emphasize on some features of the method presented here,
which are crucial for the computation of secondary tori:
\begin{itemize}
\item The method does not require neither the system to be close to integrable nor to
be written in action-angle variables.

\item The modification of the invariance equations \eqref{invariancemaps} and \eqref{invariancevf} allows to adjust some global
averages required to solve the Newton equations (see \cite{FontichLS09b}).

\item The periodicity of the function $\widetilde K$ can be adjusted by the matrix $I$ introduced in \eqref{lifting}. Hence, the rank of the matrix $I$ has to be
chosen according to the number of contractible directions.

\end{itemize}

\subsection{Equations for the invariant whiskers}\label{sec:iewhiskers}

Invariant tori with $\ell <d$ may have associated invariant bundles and
whiskers. We are interested in computing the invariant manifolds
which contain the torus and are tangent to the invariant bundles of
the linearization around the torus.  This includes the stable and
unstable manifolds but also invariant manifolds associated to
other invariant bundles of the linearization, such as the slow
manifolds, associated to the less contracting directions.

Using the parameterization method, it is natural to develop
algorithms for invariant manifolds tangent to invariant sub-bundles
that satisfy a non-resonance condition (see \cite{CabreFL03a}). This
includes as particular cases, the stable/unstable manifolds, the
strong stable and strong unstable ones as well as some other slow
manifolds satisfying some non-resonance conditions.

To avoid lengthening the paper, we restrict in this paper just to the one-dimensional manifolds (see
Section \ref{sec:rank1im}), where we do not need to deal with resonances as it is the case in higher dimensions. We think that, considering this
particular case, we can state in a more clear and simpler way the
main idea behind the algorithms. We will come back to the study of
higher dimensional manifolds in future work. 

\subsubsection{Invariant manifolds of rank 1}\label{sec:imrank1}

We once again use a parameterization to describe the whiskers. This
amounts to finding a solution $u$ of the equations of motion under
the form
$$u^{(n)}=W(\omega n, \mu^n s)$$
in the discrete time case and
$$u(t)=W(\omega t, se^{\mu t})$$
in the continuous time case, where $ W: \torus^{\ell} \times ( V
\subset \real^{d-\ell}) \rightarrow \mathcal M $ and $\mu \in
\real$. The function $W$ has then to satisfy the following
invariance equations
\begin{equation}\label{invariance-whiskers}
\begin{split}
& F (W (\theta, s)) = W (\theta + \omega, \mu
s),  \\
\noalign{\vskip6pt} &  \partial_\omega W(\theta,s) + \mu s  \frac{\partial}{\partial s}
W(\theta,s) = (X \circ W)(\theta,s),
\end{split}
\end{equation}
for the case of maps and flows, respectively.

Note that equations \eqref{invariance-whiskers} imply that in
variables $(\theta, s)$ the motion on the torus consists of a rigid
rotation of frequency $\omega$ whereas the motion on the whiskers
consists of a contraction (or an expansion) by a constant $\mu$ ($e^{\mu}$ in
the case of flows).  We call contractive the
situation  $|\mu|<1$ for maps (or  $\mu <0$
for flows). We call expansive the case when
 $|\mu|>1$ for maps (or  $\mu >0$ for
flows). Note that if $W(\theta,s)$ satisfies
\eqref{invariance-whiskers} then $W(\theta,0)$ is a solution of the
invariance equations \eqref{invariancemaps} or \eqref{invariancevf}.

As in the case of invariant tori, it will be convenient
to consider the following modified invariance equations
\begin{equation}\label{invariance-whiskers-modified}
\begin{split}
& F (W (\theta, s)) = W (\theta + \omega, \mu s) + (J (K_0)^{-1}
DK_0)
\circ T_\omega \lambda ,  \\
\noalign{\vskip6pt} &  \partial_\omega W(\theta,s) + \mu s
\frac{\partial}{\partial s} W(\theta,s) = (X \circ W)(\theta,s) + J
(K_0)^{-1} (DX\circ K_0) \lambda ,
\end{split}
\end{equation}
where $K_0$ is, as before, a given approximate solution of the
equations \eqref{invariancemaps} and \eqref{invariancevf},
respectively.

\subsubsection{Uniqueness of solutions of the invariance
equation for whiskers}\label{sec:uniqueness}
The solutions of equations
\eqref{invariance-whiskers} are not unique. Indeed, if $W(\theta,
s)$ is a solution, for any $\sigma \in \torus^\ell$, $b \in \real$,
we have that $\tilde W(\theta, s) = W (\theta + \sigma, s b)$ is
also a solution. This non-uniqueness of the problem can be removed
by supplementing the invariance equation with a normalization
condition.

Some suitable normalization conditions (in the case of maps) that
make the solutions unique are
\begin{equation}
\label{normalizationwhiskers}
\begin{split}
&\int_{\torus^\ell} W(\theta,0) - \theta = 0, \\
&DF(W(\theta,0)) D_2 W(\theta, 0) = \mu D_2 W(\theta, 0), \\
&|| DW(\cdot, 0) || = \rho
\end{split}
\end{equation}
where $D_2 W$ denotes the derivative
with respect to the second argument,
$\rho > 0$ is any arbitratrily chosen number
and $\|.\|$ stands for a suitable norm.

The fact that the solutions of \eqref{invariancemaps} supplemented
by \eqref{normalizationwhiskers} are locally unique is proved in
\cite{FontichLS09b}. In this paper, we will see that these
normalizations uniquely determine the Taylor expansions (in $s$) of
the function $W$ whenever the first term $W_1(\theta)\equiv
D_2W(\theta, 0)$ is fixed, and we will present algorithms to perform
these computations.

The first equation in \eqref{normalizationwhiskers}
amounts to choosing the origin of coordinates in the parameterization of
the torus and, therefore eliminates the ambiguity corresponding to $\sigma$.
(Check how does \eqref{normalizationwhiskers} change when we choose $\sigma$).

The second  equation in \eqref{normalizationwhiskers} indicates that
$W_1(\theta)$ is chosen to be a vector in the hyperbolic direction.
We furthermore require that we have chosen the coordinate so that it
is an eigenvector of the expanding/contracting direction.

The third equation in \eqref{normalizationwhiskers} chooses
the eigenvalue. Equivalently, it fixes  the scale in the variables $s$.
Observe  that, setting
 $b$ amounts to multiplying  $W_1$ by $b$. Hence, setting the norm
of $DW$ sets the $b$.

{From} the mathematical point of view, all choices of $\rho$ are
equivalent.  Nevertheless, from the numerical point of view, it is
highly advantageous to choose $||W_1||$ so that the numerical
coefficients of the expansion (in $s$) of $W$ have norms that
neither grow nor decrease fast. This makes the computation more
immune to round off error since it becomes more important when we
add numbers of very different sizes.

\subsection{Fourier-Taylor discretization}
\label{sec:FT} One of the ingredients of algorithms to solve the
functional equations is to consider discretizations of functions one
searches for.

In this section, we introduce the discretizations we will use.
Roughly,  for periodic functions, we will use \emph{both} a Fourier
series discretization \emph{and}  a real discretization on a
grid. We will show that the Newton step can be decomposed into
substeps which require only $O(N)$ operations in either of the
representations. Of course, one can switch between both
representations using $O(N\ln(N))$ operations using FFT algorithms.
For the study of invariant manifolds, we will use Taylor series in
the real variables.

\subsubsection{Fourier series discretization}

Since we are seeking functions $K$ which are
periodic in the angle variable $\theta$, it is natural to discretize them
 retaining a finite number
of their Fourier coefficients
\begin{equation}\label{eq:Fourierseries}
K(\theta)=\sum_{k \in \zed^{\ell}, k \in \mathcal{O}_N} c_k e^{2 i\pi k\cdot
\theta},
\end{equation}
where
$$\mathcal{O}_N= \left \{ k \in \zed^{\ell}\,\,|\,\,|k| \leq N \right \}. $$
Since we will deal with real-valued functions, we have $c_k=\bar c_{-k}$ and
one can just consider the cosine and sine Fourier series,
\begin{equation}\label{eq:fouriersincos}
K(\theta)=a_0+\sum_{k \in \mathcal{O}_N} a_k \cos (2 \pi k \cdot
\theta) + b_k \sin (2 \pi k \cdot \theta).
\end{equation}

These Fourier discretizations have a very long history going back to
classical astronomy, but have become much more widely used with
computers and go under different name such as \emph{``automatic
differentiation''}. The manipulation of these polynomials are
reviewed in \cite{Knuth97}.  A recent review of their applications
in dynamics  -- including implementation issues and examples -- is
\cite{Haro08}.

The main shortcoming of Fourier series discretization of a function
is that they are not adaptative and that for discontinuous
functions, they converge very slowly and not uniformly. These
shortcomings are however not very serious for our applications. 

Since the tori are invariant under rigid rotations, they tend to be
very homogeneous, so that adaptativity is not a great advantage.
Also, it is known \cite{FontichLS09b} that if tori are $C^r$ for
sufficiently large $r$, they are in fact analytic.

The fact that the Fourier series converge slowly for functions with
discontinuities is a slight problem if one wants to compute tori
close to the breakdown of analyticity,  when the tori transform into
Aubry-Mather objects. Of course, when they are far from breakdown --
as it happens in many interesting problems in celestial mechanics --
the Fourier coeffients converge very fast. To perform calculations
close to breakdown, the {\sl a posteriori} theorems in
\cite{FontichLS09b} prove invaluable help to have confidence in the
computed objects.

\subsubsection{Fourier vs grid representation}

Another representation of the function $K$ is to store the values in
a regularly space grid. For functions of $\ell$ variables, we see
that if we want to use $N$ variables, we can store either the
Fourier coefficients of index up to $O( N^{1/\ell})$ or the values on
a grid of step $O(N^{-1/\ell})$.

Some operations are very fast on the real space variables, for
example multiplication of functions (it suffices to multiply values
at the points of the grid). Also,  the evaluation of $F\circ K$ is
very fast if we discretize on a grid (we just need to evaluate the
function $F$ for each of the points on the grid). Other operations
are fast in Fourier representation. For example, it is fast to shift
the functions, to take derivatives and, as we will see in
\eqref{cohomology}, to solve cohomology equations. Hence, our
iterative step will consist in the application of several
operations, all of which being fast -- $O(N)$ -- either in Fourier mode
representation or in a grid representation. Of course, using the
Fast Fourier Transform, we can pass from a grid representation to
Fourier coefficients or viceversa in $O(N\ln N)$ operations. There
are extremely efficient implementations of the FFT algorithm that
take into account not only operation counts but also several other
characteristics (memory access, cache, etc.) of modern computers.

\subsubsection{Fourier-Taylor series}

For the computation of whiskers of invariant tori, we will use
Fourier-Taylor expansions of the form
\begin{equation}\label{eq:FourierTaylor}
W(\theta,s)= \sum_{n=0}^{\infty} W_n(\theta) s^n,
\end{equation}
where $W_n$ are $1$-periodic functions in $\theta$ which we will
approximate using Fourier series \eqref{eq:Fourierseries}.

To manipulate this type of series we will use the so
called \emph{automatic differentiation algorithms} (see
\cite{Knuth97},\cite{Haro08}). For the basic algebraic operations and the
elementary transcendental functions (exp, sin, cos, log, power,
etc.), they provide an expression for the Taylor coefficients of the
result in terms of the coefficients of each of the terms.

\subsection{Cohomology equations and Fourier discretization}
\label{cohomology}
In  the Newton step  to construct KAM tori, one faces
solving cohomology equations, that is, given a periodic
(on $\torus^\ell$) function $\eta$, we want to find another periodic function
$\varphi$  solving (the first equation is a small divisor equation
for flows and the second one for maps)
\begin{equation}\label{coboundary}
\begin{split}
&\partial_\omega \varphi = \eta,\\
\noalign{\vskip6pt} &\varphi - \varphi \circ T_\omega = \eta.
\end{split}
\end{equation}

As it is well known, equations \eqref{coboundary} have a solution
provided that
\begin{equation}\label{eq:compatibility}
\hat \eta_0 \equiv \int_{\torus^\ell} \eta =0,
\end{equation}
and that $\omega$ is Diophantine in the appropriate sense.   The
Fourier coefficients $\hat \varphi_k$ of the solution $\varphi$ of
\eqref{solution} are then given respectevely by
\begin{equation}\label{solution}
\begin{split}
\hat\varphi_k &= \frac{\hat\eta_k}{2\pi i\omega \cdot k},\\
\noalign{\vskip6pt} \hat\varphi_k &= \frac{\hat\eta_k}{1-e^{2\pi ik
\cdot \omega}},
\end{split}
\end{equation}
where $\hat \eta_k$ are the Fourier coefficients of the function
$\eta$.

Notice that the solution $\varphi$ is unique up to the
addition of a constant (the average $\hat \varphi_0$ of $\varphi$ is
arbitrary).

Equations \eqref{coboundary} and their solutions
\eqref{solution} are very standard in KAM theory (see the exposition
in \cite{Llave01}). Very detailed estimates can be found in
\cite{Russmann75}, when $\omega$ is Diophantine (which is our case).

\section{Fast Newton methods for
(possibly) whiskered  tori}\label{sec:whiskered}

In this section we develop an
efficient Newton method to solve the invariance equations
\eqref{invariancemaps}-\eqref{invariancevf} and \eqref{invariance-modifiedmaps}-\eqref{invariance-modifiedvf}. We mainly focus on the case of maps (the case for
vector fields being similar is described in the appendices).

We emphasize that the algorithm applies both to whiskered tori and
to Lagrangian tori. Indeed, the case of Lagrangian tori is
simpler. The hyperbolic part of the Lagrangian tori is just empty
so that we do not need to compute the splittings. We refer to
Algorithm~\ref{alg:combination} and Remark~\ref{rem:lagrangian}.

We will assume that the motion on the torus is a rigid rotation with
a Diophantine frequency $\omega \in \real^{\ell}$. As we have
already shown, the invariance implies that the vectors in the range
of $DK$ are invariant under $DF$. The preservation of the symplectic
structure implies that the vectors in the range of $(J \circ
K)^{-1} DK$ grow at most polynomially under iteration. We note also
that tori with an irrational rotation are co-isotropic, $(DK)^{T} (J
\circ K)^{-1} DK =0$, i.e.
\begin{equation}\label{eq:coisotropy}
\mathrm{Range} \, DK \cap \mathrm{Range} \, (J \circ K)^{-1} DK = \{
0 \},
\end{equation}
and therefore $\dim (\mathrm{Range} \, DK \oplus \mathrm{Range} \, (J
\circ K)^{-1} DK) = 2 \ell$. Hence, at any point of the invariant
torus of dimension $\ell$ with motion conjugate to a rotation, we
can find a $2 \ell$-dimensional space of vectors that grow at most
polynomially under iteration. As it is shown in \cite{FontichLS09a},
approximately invariant tori are approximately co-isotropic and the
transversality \eqref{eq:coisotropy} also holds.

The tori that we will consider are as hyperbolic as possible, given
the previous argument. That is, we will assume that there exist
directions that contract exponentially in the past or in the future,
which span the complement of the tangent to the torus and its
symplectic conjugate.

We will consider tori that have a hyperbolic
splitting
\begin{equation}\label{splitting}
T_{K(\theta)} \M = \E^c_{K(\theta)} \oplus \E^s_{K(\theta)} \oplus
\E^u_{K(\theta)},
\end{equation}
such that there exist $0<\mu_1, \mu_2 <1$, $\mu_3>1$ satisfying
$\mu_1 \mu_3 <1$, $\mu_2 \mu_3 <1$ and $C>0$ such that for all $n
\geq 1$ and $\theta \in \torus^{\ell}$

\begin{equation}\label{contraction-rates}
\begin{split}
\noalign{\vskip6pt} & v \in \E^s_{K(\theta)} \Longleftrightarrow
 | \M(n,\theta) v | \leq C \mu_1^n |v|  \qquad \forall n \geq 1\\
\noalign{\vskip6pt} & v \in \E^u_{K(\theta)} \Longleftrightarrow
| \M(n,\theta) v | \leq C \mu_2^n |v|  \qquad \forall n \leq 1 \\
\noalign{\vskip6pt}& v \in \E^c_{K(\theta)} \Longleftrightarrow |
\M(n,\theta) v | \leq C \mu_3^n |v|  \qquad \forall n \in \mathbb{Z}
\end{split}
\end{equation}
where $\M(n,\theta)$ is the cocycle with generator $Z(\theta)=DF
(K(\theta))$ and frequency $\omega$, i.e. $\M :\zed\times\torus^\ell \to GL (2d,\real)$ is given by
\begin{equation}
\M (n,\theta) =
\begin{cases}
Z (\theta  + (n-1)\omega) \cdots Z(\theta)&n\ge 1,\\
\Id &n=0,\\
Z^{-1} (\theta + (n+1)\omega) \cdots Z^{-1}(\theta)&n\le 1.
\end{cases}
\label{cocycle1}
\end{equation}

We will also assume that 
\begin{equation}\label{iswhiskered}
\dim \E^{c}_{K(\theta)}= 2 \ell,  \qquad \dim
\E^{s}_{K(\theta)}= \dim \E^{u}_{K(\theta)}= d - \ell.
\end{equation} 

The  assumption~\eqref{iswhiskered} 
implies that the only non-hyperbolic directions are
those spanned by the tangent to the torus and its symplectic
conjugate, that is,  there are no elliptic directions
except those that are forced by the symplectic structure
and the fact that the motion on the torus is a rotation.

We associate to the splitting \eqref{splitting} the projections
$\Pi^c_{K(\theta)}$, $\Pi^s_{K(\theta)}$ and $\Pi^u_{K(\theta)}$
over the invariant spaces $\E^c_{K(\theta)}$, $\E^s_{K(\theta)}$ and
$\E^u_{K(\theta)}$.

It is important to note that since $F$ is symplectic (i.e. $F^* \Omega=\Omega$), for all $n \geq 1$ and $n  \leq -1$
\[ \Omega (u,v)= \Omega (DF^{n} u, DF^n v),\]
so that, if $u,v$ have rates of decrease, by taking limits in the
appropriate direction we obtain that $\Omega$ is zero. That is, we
get
\begin{eqnarray*}
\Omega (\E^s, \E^s)=0, &  & \Omega (\E^u, \E^u)=0, \\
\Omega (\E^c, \E^s)=0, &  & \Omega (\E^c, \E^u)=0.
\end{eqnarray*}
Therefore, we have
\begin{eqnarray*}
(J  (K(\theta)))^{-1}  \E^c_{K(\theta)} & = & \E^c_{K(\theta)}, \\
(J  (K(\theta)))^{-1} \E^s_{K(\theta)} & = & \E^u_{K(\theta)}, \\
(J  (K(\theta)))^{-1} \E^u_{K(\theta)} & = & \E^s_{K(\theta)}.
\end{eqnarray*}


In \cite{HuguetLS10c}, we provide a method to compute the
rank-1 bundles by iterating the cocycle. Of course, once we have
computed the vector spanning the rank-1 (un)stable
bundle it is very easy to obtain the projections. In Section
\ref{sec:projections} we discuss an alternative to compute the
projections by means of a Newton method. In this case we do not need
to assume that the bundle is 1-dimensional.

\subsection{General strategy of the Newton method to solve the invariance
equation}\label{generalst}

In this section we will design a  Newton method to solve the
invariance equation \eqref{invariancemaps} and the modified one
\eqref{invariance-modifiedmaps}, and discuss several algorithms to
deal with the linearized equations.

We first define the following concept of approximate solution.
\begin{defn}
We say that $K$ (resp. $(K,\lambda)$) is an approximate solution of
equation \eqref{invariancemaps} (resp.
\eqref{invariance-modifiedmaps}) if
\begin{equation}\label{invarianceapprox}
\begin{split}
& F\circ K - K\circ T_\omega = E, \\
\noalign{\vskip6pt} & ( \mathrm{resp.} \quad F\circ  K - K\circ
T_\omega - ((J \circ K_0)^{-1} DK_0) \circ T_\omega \lambda = E ),
\end{split}
\end{equation}
where $E$ is small.
\end{defn}

The Newton method consists in computing $\Delta$ in such a way that
setting ${K\gets K+\Delta}$ and expanding the LHS of
\eqref{invarianceapprox} in $\Delta$ up to order $\|\Delta\|^2$, it
cancels the error term $E$.

\begin{remark}
Throughout the paper, we are going to denote $\|.\|$ some norms in
functional spaces without specifying what they are exactly.
We refer the reader to  \cite{LlaveGJV05,FontichLS09b}, where the
whole theory is developed and the convergence of the algorithms is
proved. Recall that one of the key ideas of KAM theory is that the
norms are modified at each step.
\end{remark}

Performing a straightforward calculation, we obtain that the Newton
procedure to solve equation \eqref{invariancemaps} and
\eqref{invariancevf}, given an approximate solution $K$, consists in
finding $\Delta$ satisfying
\begin{equation}\label{Newton}
(DF\circ K)\Delta - \Delta \circ T_\omega =  - E.
\end{equation}
For the modified invariance equation
\eqref{invariance-modifiedmaps}, given an approximate solution
$(K,\lambda)$, the Newton method consists in looking for
$(\Delta,\delta)$ in such a way that $K+\Delta$ and $\lambda
+\delta$ eliminate the error in first order. The linearized equation
in this case is
\begin{equation}\label{Newton-modified}
(DF\circ K)\Delta  - \Delta \circ T_\omega - ((J \circ K_0)^{-1}
DK_0) \circ T_\omega \delta =-E,
\end{equation}
where one can take $K_0=K$.

As it is well known, the Newton method converges quadratically in
$\|E\|$ and the error $\widetilde E$ at step $K+\Delta$ is such that
$$\|\widetilde E\| \le C\|E\|^2$$
where $E$ is the error at the previous step.

In order to solve the linearized equations \eqref{Newton} and
\eqref{Newton-modified}, we will first project them on the invariant
subspaces $\E^{c}, \E^{u}$ and $\E^s$, and then solve an equation
for each subspace.

Thus, let us denote
\begin{equation}\label{desuc}
\begin{split}
\noalign{\vskip6pt} & \Delta^{s,c,u}(\theta)=\Pi^{s,c,u}_{K(\theta)} \Delta (\theta), \\
\noalign{\vskip6pt} & E^{s,c,u}(\theta)=\Pi^{s,c,u}_{K(\theta+\omega)} E (\theta), \\
\end{split}
\end{equation}
such that
$\Delta(\theta)=\Delta^s(\theta)+\Delta^c(\theta)+\Delta^u(\theta)$.
Then, by the invariant properties of the splitting, the linearized
equations for the Newton method  \eqref{Newton} and
\eqref{Newton-modified} split into:
\begin{equation}\label{projected-eq}
\begin{split}
DF (K(\theta)) \Delta^c (\theta) - \Delta^c \circ T_{\omega}
(\theta)  & = - E^c(\theta), \\
\noalign{\vskip6pt} DF (K(\theta)) \Delta^s (\theta) - \Delta^s
\circ T_{\omega}
(\theta)  & = - E^s(\theta), \\
\noalign{\vskip6pt} DF (K(\theta)) \Delta^u (\theta) - \Delta^u
\circ T_{\omega}
(\theta)  & = - E^u(\theta), \\
\end{split}
\end{equation}
and
\begin{equation}\label{projected-eq-modified}
\begin{split}
DF (K(\theta)) \Delta^c (\theta) - \Delta^c \circ T_{\omega}
(\theta) + \Pi^c_{K(\theta+\omega)} (J\circ K_0 (\theta+\omega))^{-1} DK_{0}(\theta+\omega) \delta & = - E^c(\theta) ,\\
\noalign{\vskip6pt} DF (K(\theta)) \Delta^s (\theta) - \Delta^s
\circ T_{\omega}
(\theta) + \Pi^s_{K(\theta+\omega)} (J \circ K_0 (\theta+\omega))^{-1} DK_{0}(\theta+\omega) \delta & = - E^s(\theta) ,\\
\noalign{\vskip6pt} DF (K(\theta)) \Delta^u (\theta) - \Delta^u
\circ T_{\omega}
(\theta) + \Pi^u_{K(\theta+\omega)} (J \circ K_0 (\theta+\omega))^{-1} DK_{0}(\theta+\omega) \delta & = - E^u(\theta). \\
\end{split}
\end{equation}

Notice that once $\delta$ is obtained, the equations
\eqref{projected-eq-modified} on the hyperbolic spaces reduce to
equations of the form \eqref{projected-eq}. More precisely,
\begin{equation} \label{hyperbolic2}
DF (K(\theta)) \Delta^{s,u} (\theta) - \Delta^{s,u} \circ T_{\omega}
(\theta) = - \widetilde{E}^{s,u}(\theta)
\end{equation}
where
$$ \widetilde E^{s,u}= E^{s,u}(\theta) + \Pi^{s,u}_{K(\theta+\omega)} (J\circ K_0 (\theta+\omega))^{-1} DK_{0}(\theta+\omega)
\delta.
$$

Equations \eqref{projected-eq} and \eqref{projected-eq-modified} for
the stable and unstable spaces can be solved iteratively using the
contraction properties  of the cocycles on the hyperbolic spaces
given in \eqref{contraction-rates}. Indeed, a solution for equations
\eqref{hyperbolic2} is given explicitly by
\begin{equation}\label{eq:sumst}
\Delta^s(\theta)= \widetilde E^{s} \circ T_{-  \omega} (\theta) +
\sum_{k=1}^{\infty} (DF \circ K \circ T_{-\omega}(\theta) \times
\cdots \times DF \circ K \circ T_{- k \omega} (\theta) ) (
\widetilde E^{s} \circ T_{- (k + 1) \omega} (\theta) )
\end{equation}
for the stable equation, and
\begin{equation}\label{eq:sumunst}
\Delta^u(\theta)= - \sum_{k=0}^{\infty} (DF^{-1} \circ K(\theta)
\times \cdots \times DF^{-1} \circ K \circ T_{ k  \omega} (\theta))
( \widetilde E^{u} \circ T_{ k  \omega} (\theta) )
\end{equation}
for the unstable direction. Of course, the contraction of the
cocycles guarantees the uniform convergence of these series.

The algorithms presented in Appendix~\ref{sec:cohomology} allow us
to compute the solutions $\Delta^{s,u}$ of equations
\eqref{hyperbolic2} efficiently.

In Section~\ref{fast-whisktori} we discuss how to solve equations \eqref{projected-eq} and \eqref{projected-eq-modified} in the center direction.

Hence, the  Newton step of 
the algorithm for whiskered tori that we summarize here will
be obtained by combining several algorithms.
\begin{algorithm}\label{alg:combination}
Consider given $F$, $\omega$, $K_0$ and an approximate solution $K$
(resp. $K,\lambda$), perform the following operations:
\begin{enumerate}
\item[A)] Compute the invariant
splittings and the projections associated to the cocycle $Z(\theta)=DF \circ K(\theta)$ and $\omega$ using the
algorithms described in Section \ref{sec:projections}
 (or in \cite{HuguetLS10c}).

\item[B)]  Project the linearized equation to the hyperbolic space and use the
algorithms described in Appendix \ref{sec:cohomology} to obtain
$\Delta^{s,u}$.

\item[C)]
 Project the linearized equation on the center subspace and use the
Algorithm \ref{Fast-algorithm} in Section~\ref{fast-whisktori} to obtain $\Delta^{c}$ and $\delta$.

\item[D)] Set $K + \Delta^s + \Delta^u + \Delta^c \to K$ and $\lambda + \delta \to
\lambda$ 
\end{enumerate}
\end{algorithm}

Of course, since this is a Newton step, it will have to be 
iterated  repeatedly until one reaches solutions up to
a small tolerance error. 

We will start by some remarks on the different
steps of 
Algorithm~\ref{alg:combination} and, later, we will 
provide many more details on them. 

\begin{remark}
\label{rem:lagrangian}
It is important to remark that the above Algorithm~\ref{alg:combination}
also applies to the case of Lagrangian tori. It suffices to remark that in that, case, the center space is  the whole manifold, so that there is no need to compute the splitting.  Hence, for Lagrangian tori, the steps A) and B) of Algorithm~\ref{alg:combination} are trivial and do not need any work.
\end{remark}

\begin{remark}
The main issue of the Newton method is that it needs a good initial
guess to start the iteration. Any reasonable algorithm can be used
as an input to the Newton method. Indeed, our problems have enough
structure so that one can use Lindstedt series, variational methods,
approximation by periodic orbits, frequency methods, besides the
customary continuation methods.
\end{remark}

\begin{remark}
As we have mentioned in Remark~\ref{nonunique}, the solutions of
\eqref{invariancemaps} and \eqref{invariancevf} are not unique.
Therefore, in order to avoid dealing with non-invertible matrices in
the Newton procedure, we will impose the normalization condition
\[ \int_{\mathbb{T}^{\ell}} (K(\theta) - I(\theta)) \cdot \nu_i = 0\]
where $\{ \nu_i\}_{i=1}^L$ is a basis for $\text{Range}(I)$ ($L$
being the dimension) and $I$ is a linear function introduced in
\eqref{eq:I}.
\end{remark}

%
%
%


\subsection{Fast Newton method for (whiskered) tori: the center
directions}\label{fast-whisktori}

We present here the Newton method to solve the equations on the
center subspace in the case of maps. 

For Lagrangian
tori, the hyperbolic directions are empty and the study of 
the center direction is the only component
which is needed. Hence, the algorithms discussed in this section allow to solve, in particular, 
equations \eqref{Newton} and \eqref{Newton-modified} in the case of
Lagrangian tori.
For a discussion of  the center equations 
for Hamiltonian flows, we refer the reader to Appendix~\ref{ap:a}.

The key observation is that the linearized Newton equations
\eqref{Newton} and \eqref{Newton-modified} are  closely related to
the dynamics and therefore, we can use geometric identities to find
a linear change of variables that reduces the Newton equations to
upper diagonal difference equations with constant coefficients. This
phenomenon is often called ``automatic reducibility''.

The idea is stated in the following proposition:
\begin{proposition}[Automatic reducibility, see \cite{FontichLS09a, FontichLS09b}]\label{reducibility}
Given an approximation $K$ of the invariance equation as in
\eqref{invarianceapprox}, denote
\begin{equation}\label{notations}
\begin{split}
\alpha (\theta) & = D K (\theta)\\
\noalign{\vskip6pt}
N(\theta) & = \left( [\alpha (\theta)]^\perp \alpha (\theta)\right)^{-1}\\
\noalign{\vskip6pt}
\beta(\theta) & = \alpha (\theta) N (\theta) \\
\noalign{\vskip6pt} \gamma (\theta) & = (J\circ K (\theta))^{-1}
\beta(\theta)
\end{split}
\end{equation}
and form the following matrix
\begin{equation}\label{M-red}
M(\theta)  = [\alpha (\theta) \mid \gamma (\theta)],
\end{equation}
where by $[\cdot\mid\cdot]$ we denote the $2d\times 2\ell$ matrix
obtained by  juxtaposing the two $2d\times \ell$ matrices that are in
the arguments.

Then, we have
\begin{equation}\label{eq:reduction}
\begin{split}
& (DF\circ K(\theta)) M(\theta) = M(\theta +\omega)
\begin{pmatrix} \Id & A(\theta)\\ 0&\Id\end{pmatrix} + \widehat {E}(\theta)\\
\end{split}
\end{equation}
where
\begin{equation}\label{eq:shear}
A(\theta)= \beta(\theta+\omega)^\perp [ (DF \circ K(\theta))
\gamma(\theta) -  \gamma(\theta + \omega) ],
\end{equation}
and $\Vert \widehat E \Vert \le \Vert D E \Vert $ in the case
of \eqref{Newton} or $\Vert \widehat  E \Vert \le \Vert D E
\Vert + |\lambda|$ in the case of \eqref{Newton-modified}.
\end{proposition}

\begin{remark}
If the symplectic structure is almost-complex (i.e.
$J^2=-\Id$), we have that
\[ \beta(\theta + \omega)^\perp \gamma(\theta + \omega)=0,\]
since the torus is isotropic. Then
$A(\theta)$ has a simpler expression given by
\[ A(\theta)= \beta(\theta+\omega)^\perp (DF \circ K)(\theta) \gamma (\theta). \]
\end{remark}

Once again, we omit the definition of the norms used in the bounds for $\widehat
E$. For these precisions, we refer to the paper \cite{FontichLS09b},
where the convergence of the algorithm is established.

\begin{figure}[h]
  \begin{center}
  \psfrag{vt}{$v(\theta)$} \psfrag{ut}{$u(\theta)$} \psfrag{Kt}{$K(\theta)$}
  \psfrag{vtw}{$v(\theta+\omega)$} \psfrag{utw}{$u(\theta+\omega)$} \psfrag{Ktw}{$K(\theta + \omega)$}
  \psfrag{DfKvt}{$DF (K(\theta)) v(\theta)$}
  \includegraphics[width=100mm]{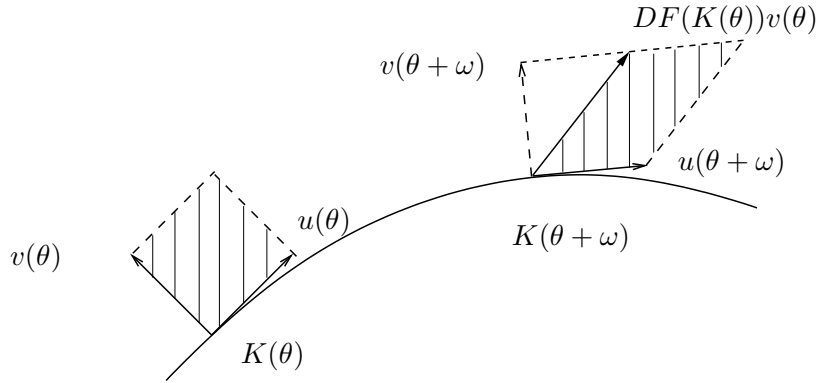}
  \caption{Geometric representation of the automatic reducibility
where $u = DK$, $v=(J \circ K)^{-1}DK N$}
  \label{fig:reduction}
  \end{center}
  \end{figure}

It is interesting to pay attention to the geometric interpretation
of the identity \eqref{eq:reduction}. Note that, taking derivatives
with respect to $\theta$ in \eqref{invarianceapprox}, we obtain that
\begin{equation*}
(DF\circ K) DK - DK \circ T_\omega =DE,
\end{equation*}
which means that the vectors $DK$ are invariant under $DF\circ K$
(up to a certain error). Moreover, $(J\circ K)^{-1} DK N$ are the
symplectic conjugate vectors of $DK$, so that the preservation of
the symplectic form clearly implies \eqref{eq:reduction}.  The
geometric interpretation of the matrix $A(\theta)$  is a shear flow
near the approximately invariant torus. See
Figure~\ref{fig:reduction}.

To be able to use the change of unknowns via the matrix
$M$ previously introduced on the center subspace, one has to ensure
that one can identify the center space $\mathcal E^c_{K(\theta)}$
with the range of $M$. This is proved in \cite{FontichLS09b}
to which we refer.

For our purposes it is
important to compute not just the invariant spaces, but also
the projections over invariant subspaces. Knowing one invariant
subspace is not enough to compute the projection, since it also depends
on the complementary space chosen.

In the following, we will see that the result stated in Proposition
\ref{reducibility} allows us to design a very efficient algorithm for
the Newton step.

Notice first that if we change the unknowns $\Delta = MW$ in
\eqref{Newton} and \eqref{Newton-modified} and we use
\eqref{eq:reduction} we obtain
\begin{equation}\label{eq:transformed}
\begin{split}
&M(\theta +\omega) \begin{pmatrix} \Id & A(\theta)\\
0&\Id\end{pmatrix}
W(\theta) - M(\theta +\omega) W(\theta +\omega) \\
\noalign{\vskip6pt} &\qquad - (J (K_0 (\theta +\omega)))^{-1} DK_0
(\theta +\omega) \delta = - E(\theta)
\end{split}
\end{equation}
Of course, the term involving $\delta$ has to be omitted when
considering \eqref{Newton}.

Multiplying \eqref{eq:transformed} by $M(\theta+\omega)^\perp
J(K(\theta+\omega))$ and using the invertibility of the matrix
$M(\theta+\omega)^\perp J(K(\theta+\omega))M(\theta+\omega)$ (see
\cite{FontichLS09a, FontichLS09b}), we are left with the system of
equations
\begin{equation}\label{eq:transformed2}
\begin{split}
& W_1 (\theta) + A(\theta) W_2(\theta) - B_1 (\theta)\delta
- W_1 (\theta +\omega) = - \widetilde E_1 (\theta)\\
\noalign{\vskip6pt} & W_2 (\theta) - W_2 (\theta +\omega) - B_2
(\theta) \delta = - \widetilde E_2 (\theta)
\end{split}
\end{equation}
where
\begin{equation*}
\begin{split}
&\widetilde E (\theta) = (M(\theta+\omega)^\perp J(K(\theta+\omega))M(\theta +\omega))^{-1} M(\theta+\omega)^\perp J(K(\theta+\omega)) E(\theta) \\
\noalign{\vskip6pt} & B(\theta)= \Big \{(M^\perp J(K)M)^{-1} M^\perp J(K)(J (K_0
))^{-1} DK_0 \Big \} \circ T_\omega (\theta)
\end{split}
\end{equation*}
and the subindices $i=1,2$ indicate symplectic coordinates.

When $K$ is close to $K_0$, we expect that $B_2$ is
close to the $\ell$-dimensional identity matrix and $B_1$ is small.

The next step is to solve equations \eqref{eq:transformed2} for $W$
(and $\delta$). Equations \eqref{eq:transformed2} are equations of the form considered in \eqref{coboundary} and
they can be solved very efficiently in Fourier space.

More precisely, the second equation of \eqref{eq:transformed2} is
uncoupled from the first one and allows us to determine $W_2$ (up to
a constant) and $\delta$. The role of the parameter $\delta$ is now
clear. It allows us to adjust some global averages that we need to
be able to solve equations \eqref{eq:transformed2}.  Indeed, we
choose $\delta$ so that the term $B_2(\theta) \delta -\widetilde
E_2$ has zero average (which is a necessary condition to solve small
divisor equations as described in Section \ref{cohomology}). This
allows us to solve equation \eqref{coboundary} for $W_2$. We then
denote
$$W_2(\theta) = \widetilde W_2 (\theta) + \overline{W}_2$$
 where
$\widetilde W_2 (\theta)$ has average zero and $\overline{W}_2 \in
\real$.

Once we have $\widetilde W_2$, we can substitute $W_2$ in the first
equation. We get $\overline{W}_2$ imposing that the average of
\[ B_1 (\theta)\delta - A(\theta) \widetilde{W}_2(\theta) - A(\theta) \overline{W}_2(\theta)  - \widetilde E_1 (\theta) \]
is zero and then we can find $W_1$ up to a constant according to
\eqref{solution}.

We therefore have the following algorithm to solve \eqref{eq:correction_split}
in the center direction,

\begin{algorithm}[ Newton step in the center direction]\label{Fast-algorithm}
Consider given $F$, $\omega$, $K_0$ and an approximate solution $K$
(resp. $K,\lambda$). Perform the following calculations
\begin{itemize}
\item[1.]
\begin{itemize}
\item[{\rm (1.1)}] Compute $F\circ K$
\item[{\rm (1.2)}] Compute $K\circ T_\omega$
\item[{\rm (1.3)}] Compute the invariant projections, $\Pi^s, \Pi^u, \Pi^c$.
\end{itemize}

\item[2.] Set $E^c = \Pi^c( F\circ K - K\circ T_\omega) $
     (resp. set $E^c= \Pi^c( F\circ K- K\circ T_\omega - (J\circ K_0)^{-1} DK_0 \lambda) $)

\item[3.] Following \eqref{notations}
\begin{itemize}
\item[(3.1)] Compute $\alpha (\theta) = DK(\theta)$
\item[(3.2)] Compute $N(\theta) = \left( [\alpha (\theta)]^\perp \alpha (\theta)\right)^{-1}$
\item[(3.3)] Compute $\beta(\theta)  = \alpha (\theta) N (\theta)$
\item[(3.4)] Compute $\gamma (\theta) = (J(K (\theta)))^{-1}\beta (\theta)$
\item[(3.5)] Compute $M(\theta) = [\alpha (\theta) \mid \gamma (\theta)]$
\item[(3.6)] Compute $M(\theta +\omega)$
\item[(3.7)] Compute $(M(\theta +\omega)^\perp J(K(\theta+\omega))M(\theta+\omega))^{-1}$
\item[(3.8)] Compute $\widetilde E (\theta) =  (M(\theta +\omega)^\perp J(K(\theta+\omega))M(\theta+\omega))^{-1}E^c(\theta)$.

We denote $\widetilde E_1, \widetilde E_2$ the components of $\widetilde E$ along
$DK$ and along $J^{-1}DK$.
\item[(3.9)] Compute
$$A(\theta) =  \beta(\theta+\omega)^\perp [ (DF \circ K(\theta))
\gamma(\theta) -  \gamma(\theta + \omega)]$$ as indicated in
\eqref{eq:shear}
\end{itemize}
\item[4.] {\rm (4.1)} Solve for $W_2$ satisfying
$$W_2 - W_2 \circ T_\omega = - \widetilde E_2 - \int_{\torus^\ell} \widetilde E_2$$
(resp.
\begin{itemize}
\item[{}]  {\rm (4.1$'$)} Solve for $\delta$ such that
$$\int_{\torus^\ell} \widetilde E_2 - \bigg[ \int_{\torus^\ell} B_2\bigg]\delta =0$$
\item[{}]  {\rm (4.2$'$)} Solve for $W_2$ satisfying
$$W_2 - W_2 \circ T_\omega = - \widetilde E_2 + B_2 \delta$$
Set $W_2$ such that the average is $0$.)
\end{itemize}
\item[5.]  {\rm (5.1)} Compute $A(\theta) W_2 (\theta)$
\item[{}]  {\rm (5.2)} Solve for $\overline{W}_2$ satisfying
$$0 = \int_{\torus^\ell} \widetilde E_1 (\theta) + \int_{\torus^\ell} A(\theta) W_2 (\theta)
+ \bigg[\int_{\torus^\ell} A(\theta)\bigg] \overline{W}_2$$
\item[{}]  {\rm (5.3)} Find $W_1$ solving
$$W_1 - W_1 \circ T_\omega = -\widetilde E_1 - A (W_2 + \overline{W}_2)$$
Normalize it so that $\int_{\torus^\ell} W_1 =0$\newline (resp.
\item[{}] {\rm (5.1$'$)} Compute $A(\theta) W_2 (\theta)$
\item[{}]  {\rm (5.2$'$)} Solve for $\overline{W}_2$ satisfying
$$0 = \int_{\torus^\ell} \widetilde E_1 (\theta)  - \int_{\torus^\ell} B_1 (\theta) \delta
+ \int_{\torus^\ell} A(\theta) W_2 (\theta) + \bigg[\int_{\torus^\ell} A(\theta)\bigg]
\overline{W}_2$$
\item[{}]  {\rm (5.3$'$)} Find $W_1$ solving
$$W_1 - W_1 \circ T_\omega = -\widetilde E_1 - A (W_2 + \overline{W}_2)
 + B_1 \delta $$
Normalize it so that $\int_{\torus^\ell} W_1 =0$.)\newline
\item[6.] The improved $K$ is $K(\theta) + M(\theta) W (\theta)$\newline
(resp. the improved $\lambda$ is $\lambda +\delta$).
\end{itemize}
\end{algorithm}

Notice that steps (1.2), (3.1), (3.6), (4.1) (resp. (4.2$'$)), (5.3)
(resp. (5.3$'$)) in Algorithm \ref{Fast-algorithm} are diagonal in
Fourier series, whereas the other steps are diagonal in the real
space representation. Note also that the algorithm only stores
vectors whose size is of order $N$.

\begin{remark} 
Using the symplectic properties of the matrix $M$, step (3.8) can be 
sped up. 

When the torus is exactly invariant we have 
that the invariant torus is co-isotropic. Hence 
$DK^\perp  J\circ K  DK = 0$. 
Hence, when the torus is invariant, we have
\[
M(\theta +\omega)^\perp J(K(\theta+\omega))M(\theta+\omega) = 
\begin{pmatrix}
0 & N \\
-N & 0 
\end{pmatrix} 
\]
so that the inverse is easy to calculate.

For the purposes of a Newton Method, we can use 
the same expression for the inverse in step (3.8) and still obtain 
a quadratically convergent algorithm. 
\end{remark}

\subsection{A Newton method to compute the
projections over invariant subspaces}\label{sec:projections}

In this section we will discuss a Newton method to compute the
projections $\Pi_{K(\theta)}^c$, $\Pi_{K(\theta)}^s$ and
$\Pi_{K(\theta)}^u$ associated to the linear spaces
$\E_{K(\theta)}^c$, $\E_{K(\theta)}^s$ and $\E_{K(\theta)}^u$ where $K$ is an (approximate) invariant torus. More precisely, we will design a  Newton method to compute $\Pi_{K(\theta)}^s$ and
$\Pi_{K(\theta)}^{cu} = \Pi_{K(\theta)}^c + \Pi_{K(\theta)}^u$.
Similar arguments allow us to design a Newton method to compute
$\Pi_{K(\theta)}^u$ and $\Pi_{K(\theta)}^{cs} = \Pi_{K(\theta)}^c +
\Pi_{K(\theta)}^s$. Then, of course, $\Pi_{K(\theta)}^c$ is given by
\begin{equation*}
\Pi_{K(\theta)}^c = \Pi_{K(\theta)}^{cs} \Pi_{K(\theta)}^{cu} =
\Pi_{K(\theta)}^{cu} \Pi_{K(\theta)}^{cs}\ .
\end{equation*}

Let us discuss first a Newton method to compute $\Pi_{K(\theta)}^s$
and $\Pi_{K(\theta)}^{cu}$. To simplify notation, from now on, we
will omit the dependence in $K(\theta)$.

Given a cocycle $Z(\theta)$ (which in our case will be $Z(\theta) =
DF(K(\theta))$), we will look for maps $\Pi^s : \torus^\ell \to
GL(2d,\real)$ and $\Pi^{cu} :\torus^\ell \to GL(2d,\real)$
satisfying the following equations:
\begin{align}
&\Pi^{cu}(\theta +\omega) Z(\theta) \Pi^s (\theta) =0,   \label{proj1}\\
&\Pi^s (\theta+\omega) Z(\theta) \Pi^{cu} (\theta) =0,   \label{proj2}\\
&\Pi^s (\theta) + \Pi^{cu} (\theta) = \Id,       \label{proj3}\\
&[\Pi^s(\theta)]^2 = \Pi^s (\theta),             \label{proj4}\\
&[\Pi^{cu}(\theta)]^2 = \Pi^{cu}(\theta),        \label{proj5}\\
&\Pi^s (\theta) \Pi^{cu} (\theta) = 0,           \label{proj6}\\
&\Pi^{cu} (\theta) \Pi^s (\theta) = 0.           \label{proj7}
\end{align}

Notice that the system of equations \eqref{proj1}--\eqref{proj7}
 is redundant. It is easy to see that equations~\eqref{proj5}, \eqref{proj6}
and \eqref{proj7} follow from equations~\eqref{proj3} and
\eqref{proj4}. Therefore, the system of equations that needs to be
solved is reduced to equations~\eqref{proj1}--\eqref{proj4}.

We are going to design a Newton method to solve
equations~\eqref{proj1}--\eqref{proj2} and use
equations~\eqref{proj3}--\eqref{proj4} as constraints. In this
context, by approximate solution of
equations~\eqref{proj1}--\eqref{proj2}, we mean a solution
$(\Pi^s,\Pi^{cu})$ such that
\begin{align}
&\Pi^{cu} (\theta +\omega) Z(\theta) \Pi^s (\theta) =
E^{cu}(\theta),
\label{proj-error1}\\
&\Pi^s (\theta +\omega) Z(\theta) \Pi^{cu} (\theta) = E^s(\theta),
\label{proj-error2}\\
&\Pi^s (\theta) + \Pi^{cu}(\theta) = \Id ,        \label{proj-error3}\\
&[\Pi^s(\theta)]^2 = \Pi^s (\theta).          \label{proj-error4}
\end{align}
where $E^{i}$ denotes the error in a certain component. Notice that the error in equation~\eqref{proj-error1} has
components only on the center and unstable ``approximated''
subspaces and we denote it by $E^{cu}$. The same happens with the
equation~\eqref{proj-error2} but on the ``approximated'' stable
subspace. We assume that $E^{cu}$ and $E^s$ are both small.

As standard in the Newton method, we will look for increments $\Delta^s$ and $\Delta^{cu}$ in such a way that
setting $\Pi^s \gets \Pi^s +\Delta^s$ and $\Pi^{cu}\gets \Pi^{cu}
+\Delta^{cu}$, the new projections solve equations~\eqref{proj1} and
\eqref{proj2} up to order $ \|E\|^2$ where $\|E\| = \|E^s\| + \|E^{cu}\|$ for some norm $\|.\|$.

The functions $\Delta^s$ and
$\Delta^{cu}$ solve the following  equations
\begin{equation} \label{Newton-proj}
\begin{split}
&\Delta^{cu} (\theta+\omega) Z(\theta) \Pi^s(\theta) +
\Pi^{cu}(\theta+\omega)
Z(\theta) \Delta^s (\theta) = - E^{cu} (\theta)\\
&\Delta^s (\theta+\omega) Z(\theta) \Pi^{cu} (\theta) +
\Pi^s(\theta+\omega) Z(\theta) \Delta^{cu} (\theta) = - E^s (\theta)
\end{split}
\end{equation}
with the constraints
\begin{align}
&\Delta^s (\theta) + \Delta^{cu} (\theta) =0 \label{constraintsa} \\
&\Pi^s (\theta) \Delta^s (\theta) + \Delta^s (\theta) \Pi^s (\theta)
= \Delta^s (\theta)\ . \label{constraintsb}
\end{align}

By equation \eqref{constraintsa} we only need to compute
$\Delta^s$ since $\Delta^{cu} = -\Delta^s$. We now work out
equations \eqref{Newton-proj}, \eqref{constraintsa} and
\eqref{constraintsb} so that we can find $\Delta^s$.

Denote
\begin{equation}\label{notproj}
\begin{split}
\Delta_s^s = \Pi^s \Delta^s,\\
\Delta_{cu}^s = \Pi^{cu} \Delta^s,
\end{split}
\end{equation}
so that
\begin{equation}\label{split-delta}
\Delta^s = \Delta_s^s + \Delta_{cu}^s.
\end{equation}
Then  equation \eqref{constraintsb} reads
\begin{equation}\label{constraints2}
\Delta_s^s (\theta) + \Delta^s (\theta) \Pi^s (\theta) = \Delta_s^s
(\theta) + \Delta_{cu}^s (\theta),
\end{equation}
or equivalently,
\begin{equation}\label{constraints-trans}
\Delta^s (\theta) \Pi^s (\theta) = \Delta_{cu}^s (\theta)\ .
\end{equation}
By \eqref{proj-error3}, \eqref{constraints-trans} and
\eqref{split-delta} we have that
\begin{equation}\label{constraints-trans2}
\Delta^s (\theta) \Pi^{cu}(\theta) = \Delta^s (\theta) - \Delta^s
(\theta) \Pi^s (\theta) = \Delta^s (\theta) - \Delta_{cu}^s (\theta)
= \Delta_s^s (\theta).
\end{equation}

Now, using \eqref{constraintsa}, equations~\eqref{Newton-proj}
transform to
\begin{equation} \label{Newton-proj-step}
\begin{split}
&-\Delta^s (\theta+\omega) Z(\theta) \Pi^s (\theta) +\Pi^{cu}
(\theta+\omega) Z(\theta) \Delta^s (\theta) = - E^{cu}(\theta),
\\
&\Delta^s (\theta+\omega) Z(\theta) \Pi^{cu} (\theta) -\Pi^s
(\theta+\omega) Z(\theta) \Delta^s (\theta) = - E^s(\theta).
\end{split}
\end{equation}

Denoting
\begin{equation*}
\begin{split}
&N_s(\theta) = \Pi^s (\theta+\omega) Z(\theta) \Pi^s (\theta),\\
&N_{cu}(\theta) = \Pi^{cu} (\theta+\omega) Z(\theta)
\Pi^{cu}(\theta),
\end{split}
\end{equation*}
and using that $\Pi^{cu}(\theta+\omega) Z(\theta) \Pi^s(\theta)$ and
$\Pi^s (\theta+\omega) Z(\theta)\Pi^{cu}(\theta)$ are small by
\eqref{proj-error1}--\eqref{proj-error2} and $\Pi^s (\theta) +
\Pi^{cu} (\theta)= \Id$ by \eqref{proj-error3}, it is enough for the
Newton method to solve for $\Delta^s$ satisfying the following
equations
\begin{equation}\label{Newton-proj-red}
\begin{split}
&-\Delta^s (\theta+\omega) \Pi^s (\theta+\omega) N_s(\theta) +
N_{cu} (\theta)\Pi^{cu} (\theta) \Delta^s (\theta) = -
E^{cu}(\theta), \\
& \Delta^s (\theta+\omega) \Pi^{cu} (\theta+\omega) N_{cu}(\theta) -
N_s (\theta)\Pi^s (\theta) \Delta^{s} (\theta) = - E^s(\theta).
\end{split}
\end{equation}
Finally, by expressions \eqref{constraints-trans} and
\eqref{constraints-trans2} and taking into account the notations
introduced in \eqref{notproj}, equations \eqref{Newton-proj-red}
read
\begin{align}
-&\Delta_{cu}^s (\theta +\omega) N_s(\theta) + N_{cu}(\theta)
\Delta_{cu}^s (\theta) = - E^{cu} (\theta),      \label{Newton3-proj1}\\
&\Delta_s^s (\theta +\omega) N_{cu}(\theta) - N_s(\theta) \Delta_s^s
(\theta) = - E^s (\theta).        \label{Newton3-proj2}
\end{align}

In Appendix~\ref{sec:cohomology}, we discussed how to solve
efficiently equations of the form
\eqref{Newton3-proj1}-\eqref{Newton3-proj2}. Notice that they are of
the form \eqref{eq:cohomology} for $A(\theta) = N_{cu}(\theta)$,
$B(\theta) = N_s (\theta)$ and $\eta(\theta) = - E^{cu}(\theta)$ in
the case of equation~\eqref{Newton3-proj1} and $A(\theta) =
N_s(\theta)$, $B(\theta) = N_{cu}(\theta)$ and $\eta (\theta) =
+E^s(\theta)$ in the case of equation~\eqref{Newton3-proj2}.
Furthermore, $\|N_s\| <1$ and $\|N_{cu}^{-1}\| <1$. Hence, they can
be solved iteratively using the fast iterative algorithms described
in Appendix~\ref{sec:cohomology}.

The explicit expressions for $\Delta_{cu}^s$ and $\Delta_s^s$ are
\begin{equation}\label{sol-proj1}
\begin{split}
\Delta_{cu}^s (\theta) &= - \Big[ N_{cu}^{-1} (\theta) E^{cu}
(\theta)
+ \sum_{n=1}^\infty N_{cu}^{-1} (\theta) \times \cdots \times \\
&\qquad N_{cu}^{-1} (\theta + n\omega) E^{cu} (\theta  +n\omega) N_s
(\theta + (n-1) \omega) \times \cdots \times N_s(\theta)\Big]
\end{split}
\end{equation}
and
\begin{equation}\label{sol-proj2}
\begin{split}
\Delta_s^s (\theta) & = E^s (\theta-\omega) N_{cu}^{-1} (\theta
-\omega)
+ \sum_{n=1}^\infty N_s (\theta -\omega) \times \cdots\times \\
&\qquad N_s (\theta -(n+1)\omega) E^s (\theta - (n+1) \omega)
N_{cu}^{-1} (\theta - (n+1)\omega ) \times\cdots\times
N_{cu}^{-1}(\theta-\omega).
\end{split}
\end{equation}

\begin{remark}
Notice that by the way $N_{cu}(\theta)$ is defined, it is a matrix
which does not have full rank. Therefore, we understand
$N^{-1}_{cu}(\theta)$ as the ``pseudoinverse'' matrix.
\end{remark}

Finally, let us check that $\Delta^s = \Delta_{cu}^s + \Delta_s^s$
also satisfies the constraints. In order to check
that constraint \eqref{constraintsb}, which is equivalent to
\eqref{constraints-trans}, is satisfied we will use the expressions
\eqref{sol-proj1} and \eqref{sol-proj2}. Notice first that
\begin{equation}\label{comp1}
N_s(\theta) \Pi^s (\theta) = N_s (\theta) \end{equation} and
\begin{equation}\label{comp2}
N_{cu}^{-1} (\theta -\omega) \Pi^s(\theta) = 0\ .
\end{equation}
Moreover, from \eqref{proj-error1} and using \eqref{proj-error4} one
can see that
\begin{equation}\label{comp3}
E^{cu} (\theta) \Pi^s (\theta) = \Pi^{cu} (\theta +\omega) Z(\theta)
[\Pi^s (\theta)]^2 = E^{cu} (\theta)\ .
\end{equation}
Then, from expressions \eqref{sol-proj1} and \eqref{sol-proj2} and
the above expression \eqref{comp1}, \eqref{comp2} and \eqref{comp3},
it is clear that
\[ \Delta^s (\theta) \Pi^s(\theta)= \Delta^s_s (\theta) \Pi^s(\theta) + \Delta^s_{cu} (\theta) \Pi^s(\theta) = 0 + \Delta^s_{cu},\]
hence, constraint \eqref{constraints-trans} is satisfied.

Now, using equation \eqref{constraintsa} we get $\Delta^{cu}
(\theta) = - (\Delta_s^s (\theta) + \Delta_{cu}^s (\theta))$ and the
improved projections are
\begin{align*}
&\widetilde\Pi^s (\theta ) = \Pi^s (\theta) + \Delta_s^s (\theta)
+ \Delta_{cu}^s (\theta)\\
&\widetilde\Pi^{cu} (\theta) = \Pi^{cu} (\theta) + \Delta^{cu}
(\theta).
\end{align*}

The new error for equations \eqref{proj1} and \eqref{proj2} is now
$\|\widetilde E\| \le C\| E\|^2$ where $\|E\| = \|E^{cu}\| +
\|E^s \|$. Of course equation \eqref{proj3} is clearly satisfied but
\eqref{proj4} is satisfied up to an error which is quadratic in
$\|E\|$. However it is easy to get an exact solution  for \eqref{proj4} and
the correction is quadratic in $\Delta^s$ (and therefore in
$\Delta^{cu}$). To do so, we just take the SVD decomposition of
$\widetilde\Pi^s$ and we set the values in the singular value
decomposition to be either 1 or 0.

In this way we obtain new projections $\Pi_{\text{new}}^s$ and
$\Pi_{\text{new}}^{cu} = \Id - \Pi_{\text{new}}^s$ satisfying
\begin{gather*}
\|\Pi_{\text{new}}^s - \widetilde \Pi^s\| < \|\Delta^s \|^2\\
\|\Pi_{\text{new}}^{cu} -\widetilde\Pi^{cu}\| < \| \Delta^{cu} \|^2,
\end{gather*}
so that the error for equations \eqref{proj1} and \eqref{proj2} is
still quadratic in $\|E\|$. Moreover, they satisfy equations
\eqref{proj4} and, of course, \eqref{proj3} exactly.

Hence, setting $\Pi^s \gets \Pi_{\text{new}}^s$ and $\Pi^{cu} \gets
\Pi_{\text{new}}^{cu}$ we can repeat the procedure described in this
section and perform another Newton step.

Consequently, the algorithm of the Newton method to compute the projections
is:
\begin{algorithm}[Computation of the projections by a Newton method]\label{alg:Newton-proj}
Consider given $F,K,\omega$ and an approximate solution
$(\Pi^s,\Pi^{cu})$ of equations \eqref{proj1}-\eqref{proj2}. Perform
the following calculations:
\begin{itemize}
\item[1.] Compute $Z(\theta) = DF\circ K(\theta)$
\item[2.] {\rm (2.1)} Compute $E^{cu}(\theta) = \Pi^{cu} (\theta+\omega)
Z(\theta) \Pi^s (\theta)$
\item[{}] {\rm (2.2)} Compute $E^s (\theta) = \Pi^s (\theta+ \omega)
Z(\theta)  \Pi^{cu}(\theta)$
\item[3.] {\rm (3.1)} Compute $N_s (\theta) = \Pi^s (\theta +\omega) Z(\theta)
\Pi^s (\theta)$
\item[{}] {\rm (3.2)} Compute $N_{cu}(\theta) = \Pi^{cu} (\theta+\omega)
Z(\theta) \Pi^{cu} (\theta)$
\item[4.] {\rm (4.1)} Solve for $\Delta_s^s$ satisfying
$$N_s (\theta) \Delta_s^s(\theta) - \Delta_s^s (\theta+\omega) N_{cu} (\theta)
= E^s (\theta)$$
\item[{}] {\rm (4.2)} Solve for $\Delta_{cu}^s$ satisfying
$$N_{cu}(\theta) \Delta_{cu}^s (\theta) - \Delta_{cu}^s (\theta +\omega)
N_s(\theta) = - E^{cu}(\theta)$$
\item[5.] {(\rm{5.1})}Compute $\widetilde\Pi^s (\theta ) = \Pi^s (\theta) + \Delta_s^s (\theta)
+ \Delta_{cu}^s (\theta)$.
\item[{}] {(\rm 5.2)} Compute the SVD decomposition of $\widetilde\Pi^s (\theta
)$: $\widetilde\Pi^s (\theta ) = U(\theta) \Sigma(\theta)
V^\perp(\theta)$.
\item[{}] {(\rm 5.3)} Set the values in $\Sigma(\theta)$ equal to the
closer integer (which will be either 0 or  1).
\item[{}] {(\rm 5.4)} Recompute $\bar \Pi^s (\theta ) = U(\theta) \Sigma(\theta)
V^\perp(\theta)$.
\item[6.] Set \enspace $\bar \Pi^s  \to \Pi^s$
\begin{itemize}
\item[{}] $\Id - \bar \Pi^s \to \Pi^{cu}$
\end{itemize}
and iterate the procedure.
\end{itemize}
\end{algorithm}

Notice that the matrix multiplication is diagonal in real space
representation, whereas the phase shift is diagonal in Fourier
space. A discussion on how to perform step 4 efficiently is given in
Appendix~\ref{sec:cohomology}.

\section{Computation of rank-1 whiskers of an invariant torus}\label{sec:rank1im}

In this section, we present algorithms to compute the whiskers
associated to an invariant torus, that is the invariant manifolds
that contain the torus and are tangent to the invariant bundles.

For the sake of simplicity and in order to state in a clear way the
main idea behind the methods we will only discuss the case when the
invariant whiskers are one-dimensional (i.e. $d-\ell=1$). The same
idea can be extended to compute invariant manifolds of any rank.
However, there are several new phenomena (resonances) that can
appear and need to be discussed. We plan to come back to this issue
in the future.

As we already mentioned in Section~\ref{sec:imrank1}, we will look for the whiskers by finding a
parameterization for them, so we will search for a  function $ W:
\torus^{\ell} \times ( V \subset \real^{d-\ell}) \rightarrow
\M $ and a scalar $\mu$ satisfying equation
\eqref{invariance-whiskers}.

We will consider three different methods to solve equation
\eqref{invariance-whiskers}. We will first discuss the order by
order method.
The other two methods are based on
the philosophy of quasi-Newton methods. Using the phenomenon of
``automatic reducibility'', we are able to design an efficient
Newton method. The first method allows to compute simultaneously the
invariant tori and the whiskers, whereas the second one assumes that
the invariant tori and the tangent bundles are already known.

We detail only the case of maps because the same ideas work for
the case of vector fields and refer the reader to the appendices for the case of flows.

Similar algorithms were developed and implemented in
\cite{HaroL06b,HaroL07} for the slightly simpler case of
quasi-periodic maps.

\subsection{The order by order method}\label{sec:obo-maps}

In this  section we adapt  the parameterization method introduced in
\cite{CabreFL05}. The convergence of the Fourier-Taylor series in this
paper 
can be easily adapted to the present case. 
 We focus on the case of maps and refer the
reader to Appendix~\ref{ap:b} for the case of flows.

We will find a solution $(W,\mu)$ of the invariance equation
\eqref{invariance-whiskers} discretizing it in Fourier-Taylor
series. Hence, we will look for $W$ as a power series
\begin{equation}\label{eq:taylor}
W(\theta,s)=\sum_{n=0}^{\infty} W_{n}(\theta) s^{n},
\end{equation}
and match similar coefficients in $s^{n}$ on both sides of equation
\eqref{invariance-whiskers}.

For $n=0$, we obtain
\begin{equation}\label{whiskers-n0}
F (W_0 (\theta)) = W_0 (\theta + \omega),
\end{equation}
which is equation \eqref{invariancemaps} for the invariant torus.
Therefore, we have $W_0(\theta)=K(\theta)$, where $K$ is a
parametrization of the invariant torus.

For $n=1$, we obtain
\begin{equation}\label{whiskers-n1}
DF \circ K (\theta) W_1(\theta) = W_1 (\theta + \omega) \mu,
\end{equation}
so that $W_1(\theta)$ is an eigenfunction with
eigenvalue $\mu$ of the operator $\mathcal{M}(1,\theta)$
defined in equation \eqref{cocycle1}.

Equation \eqref{whiskers-n1} states that the bundle spanned by $W_1$
is invariant for the linearization of $F$, and the dynamics on it is
reduced to a contraction/expansion by a constant $\mu$.
Therefore, one can compute $W_1$ and $\mu$ using the algorithms
given in Section~\ref{sec:projections}.

\begin{remark}\label{rem:bscale}
Notice that if $W_1(\theta)$ is a solution of equation
\eqref{whiskers-n1}, then $b W_1(\theta)$, for any $b \in \real$, is
also a solution. See Section \ref{sec:uniqueness} for a discussion
on how to choose $b$.
\end{remark}

For $n \geq 2$, we  obtain
\begin{equation}\label{whiskers-n2}
DF \circ K (\theta) W_n(\theta) = W_n (\theta + \omega) \mu^n +
R_n[W_0, \ldots, W_{n-1}](\theta),
\end{equation}
where $R_n$ is an explicit polynomial in $W_0,\ldots,W_{n-1}$ whose
coefficients are derivatives of $F$ evaluated at $W_0$.

Equation \eqref{whiskers-n2} can be solved provided that $\mu$ is such that $\mu^{n}$ is not in
the spectrum of the operator $\mathcal{M}(1,\theta)$. This
condition is clearly satisfied in the case of (un)stable bundles
which are one-dimensional but it can also be satisfied by other
bundles.

Equation \eqref{whiskers-n2} can be solved using the large matrix
method. It consists on considering a discretization of equation
\eqref{whiskers-n2} using Fourier series and reduce the problem to
solving a linear system in Fourier space, where the unknowns are the
Fourier coefficients of the matrix $W_n$.

There are also efficient algorithms 
which are variants of the methods devoted in the previous 
sections. 
The equation \eqref{whiskers-n2} is 
equivalent to 
\[
W_n (\theta) = (DF\circ K(\theta))^{-1}\big[  \mu^n W_n(\theta + \omega)
R_n[W_0, \ldots, W_{n-1}](\theta)\big], 
\]
which, for large enough $n$ is a contraction, so that we 
can apply the fast methods of Section~\ref{sec:cohomiter}. 
In particular Algorithm ~\ref{cohomology-iterative}.
In the case that the stable and unstable 
directions are one dimensional -- which is the one 
we discuss in this paper --  this is enough 
(remember that we always have $n \ge 2$. 
When the bundles are higher dimensional,
we may need to find a splitting corresponding to 
the cocycle generated by 
$Z(\theta) =  (DF\circ K(\theta))^{-1} \mu^n $.

\begin{remark} 
Notice that once $W_0(\theta)$ and $W_1(\theta)$ are fixed, the
solution $W_n(\theta)$ for $n \geq 2$ of equation
\eqref{whiskers-n2} is uniquely determined. It is then clear that
any analytic solution is unique. The existence of analytic solutions
is discussed in \cite{CabreFL05}.
\end{remark}

\begin{remark}
Notice that the equations to compute the new term $W_n$ do not
involve small divisors.
\end{remark}

\begin{remark}\label{rem:counterterm}
In this case, we have not considered the modified equation
\eqref{invariance-whiskers-modified} with the counterterm, because
for this method there are no obstructions to deal with as it is in
the case of the Newton method, see Remark \ref{rem:sicounterterm}.
Indeed, the vanishing Lemma in \cite{FontichLS09b} guarantees
that for exact symplectic maps $\lambda=0$ in
\eqref{invariance-whiskers-modified} once the Newton method has
converged.
\end{remark}

\subsection{A Newton method to compute simultaneously the
invariant torus and the whiskers} \label{Newton-whiskers}

In this Section we present an algorithm to solve equation
\eqref{invariance-whiskers} using a Newton method, instead of
solving it step by step as we discussed in the previous section. As
before, we only deal with the case of maps and refer the reader to
Appendix~\ref{ap:c} for the case of flows. We do not prove the
convergence of the algorithm here for sake of length (and purpose of
the paper) but it can be done using  techniques similar to those
developed in \cite{FontichLS09a, FontichLS09b}.

We start with an initial approximation $(W,\mu)$ (resp.
$(W,\mu,\lambda)$) satisfying
\begin{equation}\label{error-whiskers}
\begin{split}
&F(W(\theta,s)) - W(\theta +\omega,\mu s) = E(\theta,s)\\
&F(W(\theta,s)) - W(\theta +\omega,\mu s) - J(K_0(\theta +
\omega)^{-1} DK_0 (\theta + \omega) \lambda = E(\theta,s)
\end{split}
\end{equation}
and we look for an improved solution
\begin{equation*}
\begin{split}
&W \leftarrow W+\Delta\\
&\lambda \leftarrow  \lambda +\delta\\
&\mu \leftarrow \mu +\delta_\mu
\end{split}
\end{equation*}
by solving the linearized equation
\begin{equation}\label{eq:Newton-whiskers}
\begin{split}
&[DF\circ W(\theta,s)] \Delta (\theta,s) - \Delta (\theta
+\omega,\mu s)
-  s \partial_s W(\theta +\omega,\mu s) \delta\\
&\qquad - ((J\circ K_0 )^{-1} DK_0) \circ T_\omega (\theta)
\delta_\mu = - E(\theta, s).
\end{split}
\end{equation}

\begin{remark}\label{rem:sicounterterm}
As in the Newton method for invariant tori, the role of the
parameter $\lambda$ is to adjust some averages to solve the
equations for the case $n=1$. More precisely, $\delta$ will be
chosen in equation \eqref{eq:delta} so that equation
\eqref{order1-3} can be solved without any obstruction.
\end{remark}

We will try to solve equation \eqref{eq:Newton-whiskers} by
discretizing it in Fourier-Taylor series. Notice that equation
\eqref{eq:Newton-whiskers} is not diagonal when discretized in
Fourier-Taylor series because of the term $DF\circ W$. However,
there is a way to make it diagonal using the geometric identities
which are a direct generalization of those used for the automatic
reducibility.

We first give the idea of the automatic reducibility when $W$ is
such that
\begin{equation}\label{exact}
(F \circ W) (\theta,s)= W(\theta+\omega,\mu s).
\end{equation}

Taking derivatives with respect to $\theta$ and $s$, we see that
\begin{equation*}
\begin{split}
& DF\circ W  (\theta,s) D_{\theta}W(\theta,s)
= D_{\theta}W (\theta + \omega,\mu s),\\
& DF\circ W (\theta,s) \partial_s W (\theta,s) = \mu \partial_s
W(\theta +\omega,\mu s).
\end{split}
\end{equation*}

{From} the above equations, we read that the quantity $D_{\theta}W(\theta,s)$
remains invariant under $DF\circ W(\theta,s)$, whereas the vector
$\partial_s W(\theta,s)$ is multiplied by a factor $\mu$.

The vectors  $(J \circ W)^{-1} D_{\theta} W N $ and $(J \circ
W)^{-1} \partial_s W \tilde N$, where $N$ and $\tilde N$ are
normalization matrices (see \eqref{red-whiskers-elem}), are the
symplectic conjugate vectors of $D_{\theta} W$ and $\partial_s W $,
respectively. The preservation of the symplectic structure implies
that
\begin{equation*}
\begin{split}
&(DF\circ W(\theta,s)) (J( W(\theta, s)))^{-1} D_{\theta} W
(\theta,s) N (\theta,s)= \\
&\qquad\qquad (J( W (\theta + \omega, \mu s)))^{-1}
D_{\theta}W(\theta + \omega, \mu s)N (\theta + \omega, \mu
s)
+ A(\theta,s) D_\theta W(\theta + \omega , \mu s),\\
&(DF\circ W(\theta,s)) (J( W (\theta,s)))^{-1} \partial_s
W(\theta,s)
\tilde N (\theta,s)= \\
&\qquad\qquad \frac1{\mu} (J ( W (\theta + \omega, \mu
s)))^{-1}
\partial_s W (\theta + \omega, \mu s) \tilde N (\theta + \omega, \mu s)+ B(\theta,s) \partial_s
W(\theta + \omega, \mu s).
\end{split}
\end{equation*}
where $A(\theta,s)$ and $B(\theta,s)$ are some matrices, which will
be computed as before. See Proposition ~\ref{reducibility-whiskers}.

Therefore, one can see that in the basis $D_{\theta}W$,$(J \circ
W)^{-1} D_{\theta}W N$, $\partial_s W$, $(J \circ W)^{-1}
\partial_s W \tilde N$, the matrix $DF \circ W$ is upper triangular with constant coefficients on the
diagonal.

Following the proofs in \cite{FontichLS09b} for instance, one can
prove  the following proposition, which generalizes Proposition
\ref{reducibility}:

\begin{proposition}\label{reducibility-whiskers}
Denote
\begin{equation}\label{red-whiskers-elem}
\begin{split}
\alpha (\theta,s)   & = D_\theta W(\theta,s)\\
\beta(\theta,s)     & =  \partial_s W(\theta,s)\\
P (\theta,s)        & = \alpha(\theta,s) N (\theta,s)\\
Q(\theta,s)         & =  \beta(\theta,s) \tilde N (\theta,s)\\
N(\theta,s)         & = (\alpha (\theta,s)^\perp \alpha (\theta,s))^{-1}\\
\tilde N(\theta,s)  & = (\beta (\theta,s)^\perp \beta (\theta,s))^{-1}\\
\gamma (\theta,s)   & = (J\circ W(\theta, s))^{-1} P(\theta,s)\\
\eta (\theta,s)     & = (J\circ W(\theta, s))^{-1} Q(\theta,s)
\end{split}
\end{equation}
and form the following matrix
\begin{equation}\label{M-whiskers}
M(\theta,s) = [\alpha (\theta,s) \mid \gamma (\theta,s) \mid \beta
(\theta,s) \mid \eta (\theta,s)]
\end{equation}
where we denote by $[\cdot\mid\cdot\mid\cdot\mid\cdot]$ the
$2d\times 2d$ matrix obtained by juxtaposing the two $2d\times\ell$
matrices and the two $2d\times (d-\ell)$ matrices that are in the
arguments.

Then
$$(DF\circ W) (\theta,s) M(\theta,s) = M(\theta+\omega,\mu s) R (\theta,s) + O(E),
$$
\vskip6pt \noindent where
\begin{equation} \label{matrixtw-maps}
R(\theta,s)= \left(
\begin{picture}(150,50)
\put(10,20){$\begin{matrix} \Id &A(\theta,s)\\ \noalign{\vskip6pt}
0&\Id\end{matrix}$} \put(25,-30){$\bigcirc$}
\put(110,20){$\bigcirc$} \put(85,-30){$\begin{matrix}
\mu&B(\theta,s)\\ \noalign{\vskip6pt} 0 &1/\mu
\end{matrix}$} \put(0,0){\line(3,0){150}}
\put(75,40){\line(0,-3){90}}
\end{picture}
\right)
\end{equation}
with
\begin{equation*}
\begin{split}
& A (\theta,s)= P(\theta,s)^\perp [(DF \circ W)(\theta,s)
\gamma(\theta,s) - \gamma(\theta + \omega, \mu s)], \\
& B (\theta,s)= Q(\theta,s)^\perp [(DF \circ W)(\theta,s)
\eta(\theta,s) - \frac{1}{\mu} \eta (\theta + \omega, \mu s)], \\
\end{split}
\end{equation*}
and $E$ is the error in \eqref{error-whiskers}.
\end{proposition}

\begin{remark}
If the symplectic structure induces an almost-complex one (i.e.
$J^2=- \Id$), we have that
\begin{equation*}
\begin{split}
& A (\theta,s)= P(\theta,s)^\perp (DF \circ W)(\theta,s)
\gamma(\theta,s), \\
& B (\theta,s)= Q(\theta,s)^\perp (DF \circ W)(\theta,s)
\eta(\theta,s). \\
\end{split}
\end{equation*}
\end{remark}

Now if we change the unknowns $\Delta = MV$ in
\eqref{eq:Newton-whiskers} and multiply by $M^{-1} (\theta
+\omega,\mu s)$ the LHS, by
Proposition~\ref{reducibility-whiskers}, we are left with the
following system of equations
\begin{equation}\label{Newton-whiskers-mod}
R (\theta,s) V(\theta,s) - V(\theta + \omega,\mu s) -
C(\theta,s) \delta_u = - \widetilde E(\theta,s) + s H(\theta,s)
\delta,
\end{equation}
where $R(\theta,s)$ is given in \eqref{matrixtw-maps} and
\begin{equation*}
\begin{split}
& C (\theta,s) = M^{-1} (\theta+\omega,\mu s) (J (K_0 (\theta + \omega)))^{-1} DK_0 (\theta + \omega),\\
&\widetilde E(\theta,s) = M^{-1}(\theta+\omega, \mu s)E(\theta,s),\\
&H(\theta,s) = M^{-1} (\theta+\omega,\mu s) \partial_s W (\theta
+\omega,\mu s).
\end{split}
\end{equation*}

We expand the terms in \eqref{Newton-whiskers-mod} in power
series up to some order $L$ and match coefficients of the same order
on both sides of the equation. We use subindices to denote
coordinates and superindices to denote the order. Hence, for order
$s^0$ we have
\begin{align}
&V_1^0 (\theta) - V_1^0 (\theta +\omega) +A^0 (\theta) V_2^0(\theta)
- C_1^0 (\theta) \delta_\mu = - \widetilde E_1^0 (\theta),
\label{order0-1}\\
&V_2^0 (\theta) - V_2^0 (\theta+\omega) - C_2^0 (\theta) \delta_\mu
= -\widetilde E_2^0 (\theta), \label{order0-2}\\
&\mu V_3^0 (\theta) - V_3^0 (\theta+\omega) +B^0(\theta)
V_4^0(\theta) - C_3^0 (\theta) \delta_\mu = - \widetilde E_3^0
(\theta),
\label{order0-3}\\
&\frac1{\mu} V_4^0 (\theta) - V_4^0 (\theta+\omega) - C_4^0
(\theta)\delta_\mu = - \widetilde E_4^0 (\theta) \label{order0-4}.
\end{align}

Notice that \eqref{order0-1} and \eqref{order0-2} can be solved
using Algorithm \ref{Fast-algorithm}.
Hence, we determine $V_1^0$, $V_2^0$
and $\delta_\mu$. Once we know $\delta_\mu$, we can solve uniquely
for $V_3^0$ and $V_4^0$ in equations \eqref{order0-3} and
\eqref{order0-4}. These equations do not have any small divisors nor
obstructions since $|\mu| \neq 1$.

For order $s^1$ we have
\begin{align}
& V_1^1 (\theta) - \mu V_1^1 (\theta +\omega) + A^0(\theta)
V_2^1(\theta)
+ A^1 (\theta) V_2^0 (\theta),\label{order1-1}\\
&\qquad \qquad = -\widetilde E_1^1 (\theta) +\delta H_1^0(\theta)
+\delta_\mu C_1^1(\theta)
\notag\\
& V_2^1 (\theta) - \mu V_2^1 (\theta+\omega) = - \widetilde
E_2^1 (\theta) + \delta H_2^0 (\theta)
+ \delta_\mu C_2^1 (\theta), \label{order1-2}\\
&\mu V_3^1 (\theta) - \mu V_3^1 (\theta+\omega) +B^0
(\theta)
V_4^1(\theta) + B^1 (\theta) V_4^0 (\theta)\label{order1-3}\\
&\qquad\qquad =-\widetilde E_3^1 (\theta) + \delta H_3^0 (\theta) +
\delta_\mu C_3^1(\theta),
\notag\\
&\frac1{\mu} V_4^1 (\theta) -\mu V_4^1 (\theta +\omega) =
-\widetilde E_4^1 (\theta) +\delta H_4^0 (\theta) +\delta_\mu C_4^1
(\theta). \label{order1-4}
\end{align}

Notice that once we choose $\delta$, equations \eqref{order1-1} and
\eqref{order1-2} are uniquely solvable for $V_1^1$ and $V_2^1$.
Recall that $\delta_\mu$ is known because it has been computed in
the case of order~0 equations.

Similarly, equation~\eqref{order1-4} do not involve small
divisors nor obstructions. However, equation~\eqref{order1-3} does
have obstructions and small divisors. In order to overcome this
problem, we denote by $F$ and $G$ the solutions of
\begin{gather*}
\frac1{\mu} F(\theta) - \mu F(\theta +\omega) = H_4^0 (\theta),\\
\frac1{\mu} G(\theta) - \mu G(\theta +\omega) = D_4^1
(\theta)
\end{gather*}
where
$$D_4^1 (\theta) = - \widetilde E_4^1 (\theta) +\delta_\mu C_4^1(\theta).$$
This gives
$$V_4^1 (\theta) = \delta F(\theta) + G(\theta).$$
Taking the average of equation \eqref{order1-3}, we have that
\begin{equation}\label{eq:delta}
\overline{D_3^1} + \delta \overline{H_3^0} - \overline{B^0F} \delta
- \overline{B^0 G} - \overline{B^1 V_4^0} = 0,
\end{equation}
so we can solve for $\delta$ provided that $\overline{H_3^0} -
\overline{B^0F} \ne0$.

The other orders do not have any problem. For $s^n$, with $n
\geq 2$, we have
\begin{equation}\label{order-n}
\begin{split}
&V_1^n (\theta) - \mu^n V_1^n (\theta +\omega) + \sum_{k=0}^n
A^{n-k} (\theta) V_2^k (\theta) = - \widetilde E_1^n (\theta) +
\delta H_1^{n-1} (\theta) +
\delta_\mu C_1^n (\theta),\\
\noalign{\vskip6pt}
&V_2^n (\theta) - \mu^n V_2^n (\theta+\omega) = - \widetilde
E_2^n (\theta) + \delta H_2^{n-1}(\theta)
+ \delta_\mu C_2^n (\theta), \\
\noalign{\vskip6pt}
&\mu V_3^n (\theta) - \mu^n V_3^n (\theta+\omega) +
\sum_{k=0}^n B^{n-k}(\theta) V_4^k (\theta) = - \widetilde E_3^n
(\theta) + \delta H_3^{n-1}(\theta)
+ \delta_\mu C_3^n (\theta),\\
\noalign{\vskip6pt}
&\frac1{\mu} V_4^n (\theta) - \mu^n V_4^n (\theta +\omega) =
- \widetilde E_4^n (\theta) + \delta H_4^{n-1} (\theta) + \delta_\mu
C_4^n (\theta),
\end{split}
\end{equation}
and equations~\eqref{order-n} can be solved uniquely for $V_1^n$,
$V_2^n$, $V_3^n$ and $V_4^n$, for $n=2,\ldots , L$, where $L$ is the degree for the Taylor expansion. Hence, we have obtained
$\delta,\delta_\mu \in \real$ and
$$V (\theta,s) = \sum_{n=0}^L V^n (\theta) s^n,$$
so that the improved solution is
\begin{equation*}
\begin{split}
& W \gets W + MV,\\
& \lambda \gets \lambda +\delta,\\
& \mu \gets \mu +\delta_\mu .
\end{split}
\end{equation*}

\begin{remark}
For $L=1$, the algorithm allows us to compute
simultaneously the invariant torus and the associated linear
subspaces.
\end{remark}

The algorithm for the simultaneous computation of the whiskers and
the invariant torus is

\begin{algorithm}[Computation of tori and whiskers simultaneously]
\label{alg:wt-maps} Consider given $F$, $\omega$, $K_0$ and a fixed
order $L$. Given an approximate solution $(W,\mu,\mu)$, perform
the following calculations
\begin{itemize}
\item[1.] Compute $E(\theta,s) = F\circ W(\theta,s) - W(\theta +\omega,
\mu s) - ((J\circ K_0)^{-1} DK_0) (\theta +\omega) \mu $
\item[2.] Compute
\begin{itemize}
\item[(2.1)] $\alpha (\theta,s) = D_\theta W(\theta,s)$
\item[(2.2)] $\beta (\theta,s) = \partial_s W (\theta,s)$
\item[(2.3)] $N(\theta,s) = [\alpha (\theta,s)^\perp \alpha (\theta,s)]^{-1}$
\item[(2.4)] $\tilde N(\theta,s) = [\beta (\theta,s)^\perp \beta (\theta,s)]^{-1}$
\item[(2.5)] $P(\theta,s) = \alpha (\theta,s) N (\theta,s)$
\item[(2.6)] $Q(\theta,s) = \beta (\theta,s) \tilde N (\theta,s)$
\item[(2.7)] $\gamma (\theta,s) = (J\circ W (\theta,s))^{-1} P(\theta,s)$
\item[(2.8)] $\eta (\theta,s) = (J\circ W (\theta,s))^{-1} Q(\theta,s)$
\item[(2.9)] $M(\theta,s) = [\alpha (\theta,s) \mid \gamma (\theta,s) \mid
\beta (\theta,s) \mid \eta (\theta,s)]$
\item[(2.10)] $[M(\theta,s)]^{-1}$
\end{itemize}
\item[3.]  Compute
\begin{itemize}
\item[(3.1)] $ C (\theta,s) = M^{-1} (\theta+\omega,\mu s) (J (K_0 (\theta + \omega)))^{-1} DK_0 (\theta +
\omega)$
\item[(3.2)] $\widetilde E(\theta,s) = M^{-1} (\theta+\omega,\mu s)
E(\theta,s)$
\item[(3.3)] $H(\theta,s) = M^{-1} (\theta +\omega,\mu s) \alpha
(\theta +\omega,\mu s)$
\end{itemize}
\item[4.] Compute
\begin{itemize}
\item[(4.1)] $A(\theta,s) = P(\theta,s)^\perp [(DF \circ W)(\theta,s)
\gamma(\theta,s) - \gamma(\theta + \omega, \mu s)]$
\item[(4.2)] $B(\theta,s) = Q(\theta,s)^\perp [(DF \circ W)(\theta,s)
\eta(\theta,s) - \frac{1}{\mu} \eta (\theta + \omega, \mu
s)]$
\end{itemize}
\item[5.]
\begin{itemize}
\item[(5.1)] Solve for $\delta_\mu$ satisfying
$$\int_{\torus^\ell} \widetilde E_2^0 - \bigg[ \int_{\torus^\ell} C_2^0\bigg] \delta_\mu =0$$
\item[(5.2)] Solve for $V_2^0$ satisfying
$$V_2^0 - V_2^0 \circ T_\omega = -\widetilde E_2^0 + C_2^0 \delta_\mu$$
Set $V_2^0$ such that the average is $0$.
\end{itemize}
\item[6.]
\begin{itemize}
\item[(6.1)] Compute $A^0 (\theta) V_2^0 (\theta)$
\item[(6.2)] Solve for $\bar V_2^0$ satisfying
$$\int_{\torus^\ell} \widetilde E_1^0 - \int_{\torus^\ell} C_1^0 (\theta)\delta_\mu
+ \int_{\torus^\ell} A^0 V_2^0 + \bigg[ \int_{\torus^\ell} A^0\bigg] \bar V_2^0=0$$
\item[(6.3)] Set $V_2^0 = V_2^0 + \bar V_2^0$
\item[(6.4)] Solve for $ V_1^0$ satisfying
$$V_1^0 - V_1^0 \circ T_\omega = - \widetilde E_1^0 - A^0 V_2^0
+ C_1^0 \delta_\mu$$
\item[(6.5)] Normalize so that $\int_{\torus^\ell} V_1^0 =0$
\end{itemize}
\item[7.] Solve for $V_4^0$ satisfying
$$\frac1{\mu} V_4^0 - V_4^0 \circ T_\omega =
- \widetilde E_4^0 + C_4^0 \delta_\mu$$
\item[8.] Solve for $V_3^0$ satisfying
$$\mu V_3^0 - V_3^0 \circ T_\omega = - \widetilde E_3^0
+ C_3^0 \delta_\mu - B^0 V_4^0$$
\item[9.]
\begin{itemize}
\item[(9.1)] Solve for $F$ satisfying
$$\frac1{\mu} F -\mu F\circ T_\omega = H_4^0$$
\item[(9.2)] Solve for $G$ satisfying
$$\frac1{\mu} G - \mu G\circ T_\omega = -\widetilde E_4^1
+\delta_\mu C_4^1$$
\item[(9.3)] Solve for $\delta$ satisfying
$$\left( - \overline{\widetilde E_3^1} + \delta_\mu \overline{C_3^1}
- \overline{B^0 G} - \overline{B^1 V_4^0} \right) +\delta
(\overline{H_3^0} - \overline{B^0 F}) =0$$
\item[(9.4)] Set $V_4^1 = \delta F+G$
\end{itemize}
\item[10.\ \ ]
\begin{itemize}
\item[(10.1)] Solve for $V_3^1$ satisfying
$$\mu V_3^1 - \mu V_3^1 \circ T_\omega = -\widetilde E_3^1
+ \delta H_3^0 + \delta_\mu C_3^1 - B^0 V_4^1 - B^1 V_4^0$$
\item[(10.2)] Normalize so that $\int_{\torus^\ell} V_3^1=0$
\item[(10.3)] Solve for $V_2^1$ satisfying
$$V_2^1 - \mu V_2^1 \circ T_\omega = - \widetilde E_2^1 + \delta H_2^0
+\delta_\mu C_2^1$$
\item[(10.4)] Solve for $V_1^1$ satisfying
$$V_1^1 - \mu V_1^1 \circ T_\omega = -\widetilde E_1^1
+\delta H_1^0 + \delta_\mu C_1^1 - A^0 V_2^1 - A^1 V_2^0$$
\end{itemize}
\item[11.] For $n=2\ldots L$ do
\begin{itemize}
\item[(11.1)] Solve for $V_2^n$ satisfying
$$V_2^n - \mu^n V_2^n \circ T_\omega = -\widetilde E_2^n (\theta)
+ \delta H_2^{n-1} + \delta_\mu C_2^n$$
\item[(11.2)] Compute
$$\tilde A^n = \sum_{k=0}^n A^{n-k} V_2^k$$
\item[(11.3)] Solve for $V_1^n$ satisfying
$$V_1^n - \mu^n V_1^n \circ T_\omega = - \widetilde E_1^n
+ \delta H_1^{n-1} + \delta_\mu C_1^n - \tilde A^n$$
\item[(11.4)] Solve for $V_4^n$ satisfying
$$\frac1{\mu} V_4^n - \mu^n V_4^n \circ T_\omega
= - \widetilde E_4^n + \delta H_4^{n-1} + \delta_\mu C_4^n$$
\item[(11.5)] Compute
$$\tilde B^n = \sum_{k=0}^n B^{n-k} V_4^k$$
\item[(11.6)] Solve for $V_3^n$ satisfying
$$\mu V_3^n - \mu^n V_3^n \circ T_\omega
= - \widetilde E_3^n + \delta H_3^{n-1} + \delta_\mu C_3^n - \tilde
B^n$$
\end{itemize}
\item[12.] Compute
$$V(\theta,s) = \sum_{n=0}^L V^n (\theta) s^n$$
\item[13.] Set \quad $W\gets W+MV$
\item[{}] \qquad\quad $\mu \gets \mu +\delta$
\item[{}] \qquad\quad $\mu \gets \mu +\delta_\mu$
\end{itemize}
\end{algorithm}

\subsection{A Newton method to compute the whiskers}

Assuming that we have computed an invariant torus $K(\theta)$ with
the associated stable direction $V^s(\theta)$ (resp. unstable
direction $V^u(\theta)$) and the rate of contraction (resp.
expansion) $\mu$, we can use a Newton method to compute the
whiskers. We concentrate on the case of maps, referring to Appendix
\ref{ap:d} for the flows.

We consider the invariance equation~\eqref{invariance-whiskers}, and
we assume that we have an initial approximation $W$ for the whiskers, expressed as a power series
$$W(\theta,s) = \sum_{n=0}^\infty W^n (\theta) s^n$$
 and such that
$$W^0 (\theta) = K(\theta)\ \text{ and }\ W^1 (\theta) = V^s(\theta)$$
(the unstable case is analogous).

\begin{remark}
Again, as discussed in Remark \ref{rem:counterterm}, we do not
need to consider the modified invariance equation
\eqref{invariance-whiskers-modified} with the counterterm. The fact that $K(\theta)$ is a solution of equation
\eqref{invariancemaps} for exact symplectic maps, implies that
$\lambda=0$ in \eqref{invariance-whiskers-modified}. This is
guaranteed by the vanishing Lemma in \cite{FontichLS09b}.
\end{remark}

Then, it is clear that the error $E$ for the initial approximation
$W$ is such that
$$E(\theta,s) = \sum_{n\ge2} E^n (\theta) s^n,$$
since this is exact for the terms of order~0 and 1.

For a given function $G: \torus^{\ell} \times \real \rightarrow
\mathcal M$, we denote
\begin{equation}\label{eq:trunc}
G(\theta,s) = G^{[<L]} (\theta,s) + G^{[\ge L]} (\theta,s)
\end{equation}
where
\[
G^{[<L]}(\theta,s) = \sum_{n=0}^{L-1} G^n (\theta) s^n, \qquad
G^{[\ge L]} (\theta,s) = \sum_{n=L}^\infty G^n (\theta) s^n.
\]

Using this notation, the linearized equation for the Newton method
is
$$[DF\circ W (\theta,s)] \Delta^{[\ge2]} (\theta,s) - \Delta^{[\ge2]}
(\theta+\omega,\mu s) = - E^{[\ge2]} (\theta,s).$$

Similarly as we did in the previous section, we can perform the
change of coordinates given by the matrix $M(\theta,s)$ in
\eqref{M-whiskers} and reduce the problem to solving for
$V(\theta,s)$ the following equation, which is diagonal in
Fourier-Taylor series,
$$R(\theta,s) V^{[\ge 2]} (\theta,s) - V^{[\ge 2]} (\theta+\omega,\mu s)
= - \widetilde E^{[\ge 2]} (\theta,s),$$ with $R(\theta,s)$ given in
\eqref{matrixtw-maps}, $\widetilde E(\theta,s)=M(\theta + \omega,
\mu s)^{-1} E(\theta,s)$ and $\Delta=MV$.

Notice that in this case, we do not have to solve the system of
equations for order~0 and 1 and we can go straight to order $n\ge2$.
We use subindices to denote coordinates and superindices to denote
the order. Hence, for order $n \ge 2$ we need to solve the system of
equations
\begin{equation}\label{ordern-simple}
\begin{split}
&V_1^n (\theta) - \mu^n V_1^n (\theta +\omega) + \sum_{k=2}^n
A^{n-k} (\theta) V_2^k (\theta)
= - \widetilde E_1^n (\theta),\\
\noalign{\vskip6pt} &V_2^n (\theta) - \mu^n V_2^n (\theta
+\omega)
= - \widetilde E_2^n (\theta),\\
\noalign{\vskip6pt} &\mu V_3^n (\theta) - \mu^n V_3^n
(\theta +\omega)
+ \sum_{k=2}^n B^{n-k} (\theta) V_4^k (\theta) = - \widetilde E_3^n,\\
\noalign{\vskip6pt} &\frac1{\mu} V_4^n (\theta) -\mu^n V_4^n
(\theta +\omega) = - \widetilde E_4^n.
\end{split}
\end{equation}

The solution of \eqref{ordern-simple} for $n=2,3$
provides an exact solution of the invariance equation up to order~4.
That is, if we set
$$V^{[<4]} (\theta,s) = V^2 (\theta,s) + V^3 (\theta,s)$$
where $V^2$ and $V^3$ are obtained by solving
equations~\eqref{ordern-simple}, then the improved solution $\bar W$ is given by
$$ \bar W(\theta,s) = W(\theta,s) + M(\theta,s)V^{[<4]} (\theta,s),$$
where $M(\theta,s)$ was introduced in \eqref{M-whiskers}. The function $W$ approximates the solution of the invariance equations with
an error $\bar E$ such that
$$ \bar E(\theta,s) = \bar E^{[\ge 4]} (\theta,s).$$

This process can be iterated and at each step we solve the
invariance equation exactly up to an order which is the double of
the one we had for the initial approximation. Thus, if we assume
that the initial guess $W$ is such that the error in
\eqref{error-whiskers} satisfies
$$E = E^{[\ge L]},$$
then the modified linearized equation for the Newton method is such
that
$$R(\theta,s) V^{[\ge L]} (\theta,s) - V^{[\ge L]} (\theta+\omega,\mu s)
= - \widetilde E^{[\ge L]} (\theta,s),$$ with $R(\theta,s)$ given in
\eqref{matrixtw-maps}. If we solve the system of equations
\eqref{ordern-simple} for $n= L\ldots (2L-1)$, then the improved
$\bar W$ is
$$\bar W(\theta,s) = W(\theta,s) + M(\theta,s) V^{[< 2L]} (\theta,s),$$
with $M(\theta,s)$ as in \eqref{M-whiskers}, and the new error $\bar
E$ satisfies $\bar E(\theta,s)= \bar E^{[\ge 2L]} (\theta,s)$.

The algorithm in this case is:
\begin{algorithm}[Computation of whiskers]
Given $F$, $\omega$ as well as $K, V^s, \mu$ and an approximate
solution $W$ such that
\[ F( W(\theta,s)) - W(\theta+\omega,\mu s)=E^{[\ge L]}(\theta,s) \]
with $L \geq 2$ and $W(\theta,0)=K(\theta)$ and
$\partial_{s}W(\theta,0)=V^s(\theta)$. Perform the following
calculations:
\begin{itemize}
\item[1.]
Compute $E^{[\ge L]} (\theta,s)= F\circ W(\theta,s) -
W(\theta+\omega,\mu s)$
\item[2.] Compute
\begin{itemize}
\item[(2.1)] $\alpha (\theta,s) = D_\theta W(\theta,s)$
\item[(2.2)] $\beta (\theta,s) = \partial_s W (\theta,s)$
\item[(2.3)] $N(\theta,s) = [\alpha (\theta,s)^\perp \alpha (\theta,s)]^{-1}$
\item[(2.4)] $\tilde N(\theta,s) = [\beta (\theta,s)^\perp \beta (\theta,s)]^{-1}$
\item[(2.5)] $P(\theta,s) = \alpha (\theta,s) N (\theta,s)$
\item[(2.6)] $Q(\theta,s) = \beta (\theta,s) \tilde N (\theta,s)$
\item[(2.7)] $\gamma (\theta,s) = (J\circ W (\theta,s))^{-1} P(\theta,s)$
\item[(2.8)] $\eta (\theta,s) = (J\circ W (\theta,s))^{-1} Q(\theta,s)$
\item[(2.9)] $M(\theta,s) = [\alpha (\theta,s) \mid \gamma (\theta,s) \mid
\beta (\theta,s) \mid \eta (\theta,s)]$
\item[(2.10)] $[M(\theta,s)]^{-1}$
\end{itemize}
\item[3.] Compute
$$\widetilde E^{[\ge L]} (\theta,s) = M^{-1} (\theta +\omega,\mu s)
E^{[\ge L]} (\theta,s)$$
\item[4.] Compute
\begin{itemize}
\item[(4.1)] $A(\theta,s) = P(\theta,s)^\perp [(DF \circ W)(\theta,s)
\gamma(\theta,s) - \gamma(\theta + \omega, \mu s)]$
\item[(4.2)] $B(\theta,s) = Q(\theta,s)^\perp [(DF \circ W)(\theta,s)
\eta(\theta,s) - \frac{1}{\mu} \eta (\theta + \omega, \mu
s)]$
\end{itemize}
\item[5.] For $n=L\ldots 2L-1$ do
\begin{itemize}
\item[(5.1)] Solve for $V_2^n$ satisfying
$$V_2^n - \mu^n V_2^n \circ T_\omega = -\widetilde E_2^n (\theta)$$
\item[(5.2)] Compute
$$\tilde A^n = \sum_{k=L}^n A^{n-k} V_2^k$$
\item[(5.3)] Solve for $V_1^n$ satisfying
$$V_1^n - \mu^n V_1^n \circ T_\omega = -\widetilde E_1^n - \tilde A^n$$
\item[(5.4)] Solve for $V_4^n$ satisfying
$$\frac1{\mu} V_4^n - \mu^n V_4^n \circ T_\omega =
- \widetilde E_4^n$$
\item[(5.5)] Compute
$$\tilde B^n = \sum_{k=L}^n B^{n-k} V_4^k$$
\item[(5.6)] Solve for $V_3^n$ satisfying
$$\mu V_3^n - \mu^n V_3^n \circ T_\omega =
- \widetilde E_3^n - \tilde B^n$$
\end{itemize}
\item[6.] Compute
$$V(\theta,s) = \sum_{n=L}^{2L-1} V^n (\theta) s^n$$
\item[7.] Set $W\gets W+MV$
\end{itemize}
\end{algorithm}

\appendix

\def\thetheorem{\Alph{section}.\arabic{theorem}}

\section{Fast algorithms to solve difference equations
 with non constant coefficients}\label{sec:cohomology}

In this section we present fast algorithms to solve for $\Delta
(\theta)$ the cohomology equation with non constant coefficients
\begin{equation}\label{eq:cohomology}
A(\theta) \Delta (\theta) - \Delta (\theta+\omega) B(\theta) = \eta
(\theta)
\end{equation}
for given $A(\theta)$, $B(\theta)$ and $\eta (\theta)$ satisfying
either $\|A\|<1$, $\|B^{-1}\|<1$ or $\|A^{-1}\| <1$, $\|B\|<1$.

Equations of this form  appear in the Netwon step for whiskered tori
(See the informal description in Section~\ref{sec:overview}).
Equations of this form also appear in the calcultion of the
invariant splitting (see
\eqref{Newton3-proj1}-\eqref{Newton3-proj2}).

We will present two algorithms. The first one is an iterative method
with an accelerated convergence and the second one very fast (see
Section~\ref{sec:cohomiter}). The second one is only for the case of
one-dimensional bundles and it is faster (computations are
$O(N)$)(see Section~\ref{sec:cohom1d}).

\subsection{Fast iterative algorithms for the cohomology equation}
\label{sec:cohomiter}

In this section we will present a fast algorithm to solve equation
\eqref{eq:cohomology} using iterative methods. We refer the reader
to \cite{HuguetLS10c} where a similar idea is used to compute
iteration of cocycles.

We consider first the case $\|A^{-1}\| <1$ and $\|B\| <1$ or, more
generally, $\|A^{-1} (\theta) \| \cdot \|B(\theta)\| <1$ Then,
multiplying \eqref{eq:cohomology} by $A^{-1}(\theta)$ on the LHS, we
obtain
\begin{equation}\label{eq:ite1}
\Delta (\theta) = A^{-1} (\theta) \Delta (\theta +\omega) B(\theta)
+ A^{-1} (\theta) \eta(\theta).
\end{equation}
This is a contraction mapping and it is straightforward to iterate
it an obtain an algorithm that converges faster than exponentially.

Next, we compute $\Delta (\theta +\omega)$ by shifting
\eqref{eq:ite1} and substituting again in \eqref{eq:ite1}, so that
we get
\begin{equation*}
\begin{split}
\Delta (\theta) & = A^{-1} (\theta) \eta (\theta) \\
&\qquad + A^{-1}(\theta) A^{-1} (\theta +\omega) \eta (\theta
+\omega)
B(\theta)\\
&\qquad + A^{-1} (\theta) A^{-1} (\theta +\omega) \Delta
(\theta+2\omega) B(\theta +\omega) B(\theta).
\end{split}
\end{equation*}

Notice that if we define
\begin{equation*}
\bar \eta (\theta) = A^{-1}(\theta) \eta (\theta)\hskip.8truein
\end{equation*}
and
\begin{equation*}
\begin{split}
&A_1^{-1} (\theta) = A^{-1} (\theta) A^{-1}(\theta +\omega),\\
&B_1 (\theta) = B(\theta +\omega) B(\theta),\\
&\eta_1 (\theta) = \bar \eta (\theta) + A^{-1}(\theta) \bar
\eta(\theta + \omega) B(\theta),
\end{split}
\end{equation*}
we have that
$$\Delta (\theta) = \eta_1 (\theta) + A_1^{-1} (\theta) \Delta (\theta +
2\omega) B_1 (\theta)$$ which has the same structure as
\eqref{eq:ite1} and we can repeat the same scheme. This leads to an
iterative procedure to compute $A(\theta)$, converging
superexponentially in the number of iterations. Thus, define
\begin{equation*}
\begin{split}
& A_{n+1}^{-1} (\theta) = A_n^{-1}  (\theta) A_n^{-1} (\theta + 2^n\omega),\\
& B_{n+1} (\theta) = B_n (\theta + 2^n \omega) B_n (\theta),\\
& \eta_{n+1}(\theta) = \eta_n (\theta) + A_n^{-1}(\theta)
\eta_n(\theta + 2^n \omega) B_n(\theta),
\end{split}
\end{equation*}
for $n\ge0$, with
\begin{equation*}
\begin{split}
& A_0^{-1} (\theta) = A^{-1}(\theta),\\
& B_0 (\theta) = B(\theta),\\
&\eta_0(\theta) = \bar\eta (\theta).
\end{split}
\end{equation*}
Then,
\begin{equation*}
\Delta (\theta) = \eta_{n+1}(\theta) + A_{n+1}^{-1}(\theta) \Delta
(\theta + 2^{n+1}\omega) B_{n+1}(\theta),\,\,\, \forall \ n\ge 0
\end{equation*}
and
\begin{equation*}
\Delta (\theta) = \lim_{n\to +\infty} \eta_{n+1}(\theta).
\end{equation*}
The convergence of the algorithm is guaranteed by the contraction of
$A^{-1}$ and $B$. The cost of computing $2^n$ terms in the sum is
proportional to $n$ since it involves only $n$ steps of the
algorithm.

The iterative algorithm is the following:

\begin{algorithm}[Solution of difference equations with non constant coefficient] \label{cohomology-iterative}
Given $A(\theta)$, $B(\theta)$ such that $\| A^{-1}(\theta)\| \cdot
\| B(\theta) \| \le \kappa < 1$, and $\eta (\theta)$ perform the
following operations:
\begin{itemize}
\item[1.] Compute $\Delta (\theta) = A^{-1} (\theta) \eta (\theta)$
\item[2.] Compute
\begin{itemize}
\item[(2.1)] $\tilde\Delta (\theta) = A^{-1}(\theta)\Delta (\theta+\omega)
B(\theta) + \Delta (\theta)$
\item[(2.2)] $\tilde A^{-1} (\theta) = A^{-1} (\theta) A^{-1} (\theta+\omega)$
\item[(2.3)] $\tilde B(\theta) = B(\theta +\omega) B(\theta)$
\end{itemize}
\item[3.] Set \enspace $\tilde \Delta \to\Delta$
\begin{itemize}
\item[{}] $\tilde A \to A$
\item[{}] $\tilde B \to B$
\item[{}] $2\omega \to \omega$
\end{itemize}
\item[4.] Iterate points $2-3$
\end{itemize}
\end{algorithm}

The case when $\|A\|<1$ and $\|B^{-1}\| <1$ can be done similarly.
In this case, we multiply \eqref{eq:cohomology} by $B^{-1}(\theta)$
on the LHS so that we obtain
\begin{equation*}
\Delta (\theta +\omega) = A(\theta) \Delta (\theta) B^{-1}(\theta) -
\eta (\theta) B^{-1}(\theta).
\end{equation*}
Before applying the iterative scheme we shift by $-\omega$. In this
way, we have
\begin{equation*}
\Delta (\theta) = A(\theta') \Delta (\theta') B^{-1}(\theta') -\eta
(\theta') B^{-1} (\theta')
\end{equation*}
where $\theta' = T_{-\omega} \theta$.


The iterative algorithm in this case is

\begin{algorithm} \label{cohomology-iterative2}
Given $A(\theta)$, $B(\theta)$ $\|A(\theta)\| \| B^{-1}(\theta)\|
\le \kappa < 1$
 and $\eta(\theta)$, perform the
following operations:
\begin{itemize}
\item[1.] Compute $\Delta (\theta) = -\eta (\theta) B^{-1}(\theta)$
\item[2.] Compute
\begin{itemize}
\item[(2.1)] $\tilde\Delta (\theta) = A(\theta) \Delta (\theta-\omega)
B^{-1}(\theta) + \Delta (\theta)$
\item[(2.2)] $\tilde A(\theta) = A(\theta) A(\theta - \omega)$
\item[(2.3)] $\tilde B^{-1} (\theta) = B^{-1} (\theta-\omega) B^{-1}(\theta)$
\end{itemize}
\item[3.] Set
\begin{itemize}
\item[{}] $\tilde\Delta \to \Delta$
\item[{}] $\tilde A\to A$
\item[{}] $\tilde B\to B$
\item[{}] $2\omega\to\omega$
\end{itemize}
\item[4.] Iterate parts 2--3
\end{itemize}
\end{algorithm}

This algorithm gives $\Delta (\theta +\omega)$. Shifting it by
$-\omega$ we get $\Delta (\theta)$.

\subsection{Fast algorithm for the 1-D cohomology
equation with non-constant coefficients}\label{sec:cohom1d}

In this section we present an efficient algorithm for the
one-dimensional version of equation~\eqref{eq:cohomology}. It is an
adaptation of methods used  in \cite{Herman83}.

Consider the following equation,
\begin{equation}\label{eq:cohomology2}
\frac{A(\theta)}{B(\theta)} \Delta (\theta) - \Delta (\theta+\omega)
= \frac{\eta(\theta)}{B(\theta)}
\end{equation}
which is obtained from \eqref{eq:cohomology} multiplying by
$B^{-1}(\theta)$ (recall that in this case $B(\theta)$ is just a
number).

We will solve \eqref{eq:cohomology2} in two steps:

1. Find $C(\theta)$ and $\nu\in\real$ such that
\begin{equation}\label{step1}
\frac{A(\theta)}{B(\theta)} = \nu \frac{C(\theta)}{C(\theta+\omega)}
\end{equation}

2. Applying \eqref{step1} in \eqref{eq:cohomology2} and multiplying
by $C(\theta +\omega)$ we obtain
\begin{equation}\label{step2}
\nu C(\theta) \Delta(\theta) - C(\theta +\omega) \Delta (\theta
+\omega) = \tilde\eta (\theta)
\end{equation}
where $\tilde\eta(\theta) = C(\theta +\omega) B^{-1}
(\theta)\eta(\theta)$.

If we change the unknowns in \eqref{step2} by $W(\theta) = C(\theta)
\Delta(\theta)$, we are left with the equation
\begin{equation}\label{smalldiv2}
\nu W(\theta) - W(\theta +\omega) = \tilde \eta(\theta).
\end{equation}
Of course, if $|\nu|\ne1$, equation~\eqref{smalldiv2} can be solved
in Fourier space. That is, we can obtain the Fourier coefficients of
$W$  as:
$$\widehat W_k =
\frac{\skew2\widehat{\widetilde\eta}_k}{\nu - e^{2\pi ik\omega}}\
,$$ and the solution is unique. Notice that whenever $|\nu|=1$,
equation \eqref{smalldiv2} involves small divisors, which is not the
case for the iterative methods that will be discussed in the
following section.

Finally, once we have $W(\theta)$ we get
$$\Delta (\theta) = C^{-1}(\theta) W(\theta).$$

Step 1 can be achieved by taking logarithms of \eqref{step1}. Assume
that $A(\theta)/B(\theta)>0$, otherwise we change the sign. Then, we
get
$$\log A(\theta) - \log B(\theta) = \log \nu +\log C(\theta)
- \log C(\theta +\omega). $$

Taking $\log\nu$ to be the average of $\log A(\theta)-\log
B(\theta)$, the problem reduces to solve for $\log C(\theta)$ an
equation of the form \eqref{coboundary}. Then $C(\theta)$ and $\nu$
can be obtained by exponentiation. Notice that $\log C(\theta)$ is
determined up to a constant. We will fix it by assuming  that it has
average~$0$.

Hence, we have the following algorithm:
\begin{algorithm}[Solution of difference equations with non constant coefficient $(1D)$]\label{1Dcohomology}
Given $A(\theta)$, $B(\theta)$ and $\eta(\theta)$. Perform the
following instructions:
\begin{itemize}
\item[1.] {\rm (1.1)} Compute $L(\theta) = \log (A(\theta)) -\log(B(\theta))$
\item[{}] {\rm (1.2)} Compute $\overline L = \int_{\torus^\ell} L$
\item[2.] Solve for $L_C$ satisfying
$$L_C(\theta)- L_C\circ T_\omega (\theta) = L(\theta) - \overline L$$
as well as having zero average.
\item[3.] {\rm (3.1)} Compute $C(\theta) = \exp (L_C(\theta))$
\item[{}] {\rm (3.2)} Compute $\nu = \exp (\overline L)$
\item[4.] Compute $\tilde\eta (\theta ) = C(\theta +\omega) B^{-1}(\theta)
\eta (\theta)$
\item[5.] Solve for $W$ satisfying
$$\nu W(\theta) - W(\theta +\omega) = \tilde\eta (\theta)$$
\item[6.] The solution of \eqref{eq:cohomology} is $\Delta (\theta)
= C^{-1} (\theta) W(\theta)$
\end{itemize}
\end{algorithm}

\section{Fast Newton method for whiskered tori in Hamiltonian flows: the center
directions}\label{ap:a}

In this section, we provide the numerical algorithm to solve the
invariance equation \eqref{invariancevf} and the modified one
\eqref{invariance-modifiedvf} using a Newton method analogous to the
one described in Section \ref{fast-whisktori}.

The automatic reducibility can also be proved in this context (see
\cite{FontichLS09b}) and we provide here the algorithm only.

\begin{algorithm}[Newton step for flows in the center direction]\label{Fast-algorithm-flow}
Consider given $X=J(K) \nabla H$, $\omega$, $K_0$ and an
approximate solution $K$ (resp. $K,\lambda$). Perform the following
calculations
\begin{itemize}
\item[1.] {\rm (1.1)} Compute $ \partial_\omega K$.
\item[{}]
    {\rm (1.2)} Compute $ X \circ K$
    {\rm (1.3)} Compute the invariant projections $\Pi^c,\Pi^u,\Pi^s$
\item[2.] Set $E^c= \Pi^c(\partial_\omega K - X \circ K)$
    (resp. set $E^c= \Pi^c(\partial_\omega K - X \circ K - (J\circ K_0)^{-1} (DX \circ K_0) \lambda$))
\item[3.] Following \eqref{notations}
\begin{itemize}
\item[(3.1)] Compute $\alpha (\theta) = D K (\theta)$
\item[(3.2)] Compute $\beta (\theta) = J(K_0(\theta))^{-1}  \alpha (\theta)$
\item[(3.3)] Compute $\beta(\theta)  = \alpha (\theta) N (\theta)$
\item[(3.4)] Compute $\gamma (\theta) = (J(K (\theta)))^{-1}\beta (\theta)$
\item[(3.5)] Compute $M(\theta) = [\alpha (\theta) \mid \gamma (\theta)]$
\item[(3.6)] Compute $M(\theta +\omega)$
\item[(3.7)] Compute $(M(\theta +\omega)^\perp J(K(\theta+\omega))M(\theta+\omega))^{-1}$
\item[(3.8)] Compute $\widetilde E (\theta) = (M(\theta +\omega)^\perp J(K(\theta+\omega))M(\theta+\omega))^{-1} E(\theta)^c$
\item[(3.9)] Compute
$$S(\theta) = \beta^\perp(\theta) ( \Id_{2d} - \beta (\theta)
\alpha(\theta)^\perp ) (DX(K(\theta)) + DX(K(\theta))^\perp ) \beta
(\theta)$$.
\end{itemize}
\item[4.] {\rm (4.1)} Solve for $W_2$ satisfying
$$ \partial_\omega W_2   = - \widetilde E_2 - \int_{\torus^\ell} \widetilde E_2$$
(resp.
\begin{itemize}
\item[{}]  {\rm (4.1$'$)} Solve for $\delta$ satisfying
$$\int_{\torus^\ell} \widetilde E_2 - \bigg[ \int_{\torus^\ell} B_2\bigg]\delta =0$$
\item[{}]  {\rm (4.2$'$)} Solve for $W_2$ satisfying
$$ \partial_\omega W_2 = - \widetilde E_2 + B_2 \delta$$
Set $W_2$ such that its average is $0$.)
\end{itemize}
\item[5.]  {\rm (5.1)} Compute $S(\theta) W_2 (\theta)$
\item[{}]  {\rm (5.2)} Solve for $\overline{W}_2$ satisfying
$$\int_{\torus^\ell} \widetilde E_1 (\theta) + \int_{\torus^\ell} S(\theta) W_2 (\theta)
+ \bigg[\int_{\torus^\ell} S(\theta)\bigg] \overline{W}_2 = 0$$
\item[{}]  {\rm (5.3)} Find $W_1$ solving
$$\partial_\omega W_1 = -\widetilde E_1 - S (W_2 + \overline{W}_2) $$
Normalize it so that $\int_{\torus^\ell} W_1 =0$\newline (resp.
\item[{}] {\rm (5.1$'$)} Compute $S(\theta) W_2 (\theta)$
\item[{}]  {\rm (5.2$'$)} Solve for $\overline{W}_2$ satisfying
$$ \int_{\torus^\ell} \widetilde E_1 (\theta)  + \int_{\torus^\ell} B_1 (\theta) \delta
- \int_{\torus^\ell} S(\theta) W_2 (\theta) -
\bigg[\int_{\torus^\ell} S(\theta)\bigg] \overline{W}_2 =0 $$
\item[{}]  {\rm (5.3$'$)} Find $W_1$ solving
$$ \partial_\omega W_1 = -\widetilde E_1 - S ( W_2
+ \overline{W}_2) + B_1 \delta  $$ Normalize it so that $\int_{\torus^\ell} W_1
=0$).\newline
\item[6.] The improved $K$ is $K(\theta) + M(\theta) W (\theta)$\newline
(resp. the improved $\lambda$ is $\lambda +\delta$).
\end{itemize}
\end{algorithm}

Notice that steps (1.1), (3.1), (4.1) (resp. (4.2$'$)), (5.3)
(resp. (5.3$'$)) in Algorithm \ref{Fast-algorithm-flow} are diagonal
in Fourier series, whereas the other steps are diagonal in the real
space representation. The algorithm only stores vectors
which are of length of order $N$.

\section{The order by order method for whiskers in Hamiltonian
flows}\label{ap:b}

In this section we present the result analogous to the one described
in Section \ref{sec:obo-maps} to solve the invariance equation
\eqref{invariance-whiskers} for the whiskers in the case of
Hamiltonian flows.

As in Section \ref{sec:obo-maps} we look for $W$ as a power series
\begin{equation*}
W(\theta,s)=\sum_{n=0}^{\infty} W_{n}(\theta) s^{n},
\end{equation*}
and match similar coefficients in $s^{n}$ on both sides of equation
\eqref{invariance-whiskers}.

For $n=0$, one obtains
\begin{equation}\label{whiskers-n0-flow}
\partial_\omega W_0 (\theta) = (X \circ
W_0)(\theta)
\end{equation}
which admits the solution $W_0(\theta)=K(\theta)$, where $K$ is a
parametrization of the invariant torus.

For $n=1$, we obtain
\begin{equation}\label{whiskers-n1-flow}
\partial_\omega W_1(\theta) + W_1(\theta)
\mu = (DX \circ K(\theta)) W_1(\theta) ,
\end{equation}
from where we read that $W_1(\theta)$ is an eigenfunction with
eigenvalue $-\mu$ of the operator $\mathcal{L}_\omega$
\[\mathcal{L}_\omega := \partial_\omega - DX \circ K(\theta).\]

Again, we note that, multiplying a solution of
\eqref{whiskers-n1-flow} by a scalar $b \in \real$, we also obtain a
solution. See Remark \ref{rem:bscale}.

 For $n \geq 2$, we  obtain
\begin{equation}\label{whiskers-n2-flow}
\partial_\omega W_n(\theta) + W_n(\theta) n
\mu= (DX \circ K (\theta))W_n(\theta) + R_n(W_0, \ldots,
W_{n-1}),
\end{equation}
where $R_n$ is an explicit polynomial in $W_0,\ldots,W_{n-1}$ whose
coefficients are derivatives of $X$ evaluated at $W_0=K$.

Notice that, in this case, equation \eqref{whiskers-n2-flow} can be
solved provided that $n\mu$ is not in the spectrum of the
operator $\mathcal{L}_\omega$ (this is a non-resonance condition
which is clearly satisfied since the stable spaces are
1-dimensional). As in the case of maps, the previous equation can be
solved using the large matrix method.

\section{A Newton method to compute simultaneously the invariant torus and the whiskers for flows} \label{ap:c}

As in Section \ref{Newton-whiskers}, we start with an initial
approximation $(W,\mu)$ (resp. $(W,\mu,\mu)$ for the invariance
equation \eqref{invariance-whiskers} (resp.
\eqref{invariance-whiskers}), that is
\begin{equation}\label{error-whiskers-flow}
\begin{split}
&  X(W(\theta,s)) - \left ( \partial_\omega
+ \mu s
\frac{\partial}{\partial s} \right ) W(\theta,s) = E(\theta,s), \\
&  X(W(\theta,s))- \left ( \partial_\omega + \mu s
\frac{\partial}{\partial s} \right ) W(\theta,s) - (J^{-1} DX) \circ
K_0 \lambda =E(\theta,s)
\end{split}
\end{equation}
and we will look for an improved solution
\begin{equation*}
\begin{split}
&W\to W+\Delta\\
&\lambda \to  \lambda +\delta\\
&\mu \to \mu +\delta_\mu
\end{split}
\end{equation*}
by solving the linearized equation
\begin{equation}\label{eq:Newton-whiskers-flow}
\begin{split}
& (DX \circ (W(\theta,s)) \Delta (\theta,s) - \left (
\partial_\omega + \mu s \frac{\partial}{\partial s} \right )
\Delta (\theta,s) -  s
\frac{\partial}{\partial s} W (\theta,s) \delta \\
&\qquad -(J^{-1} DX) \circ K_0 \delta_{\mu} = -E(\theta,s).
\end{split}
\end{equation}

Once again, we will use a reducibility argument similar to the
automatic reducibility of Lagrangian tori. This will lead to a
diagonal equation. Applying the operators $D_\theta$ and
$\partial_s$ to equations \eqref{error-whiskers-flow}, we obtain
\begin{equation*}
\begin{split}
&
DX(W(\theta,s)) D_{\theta} W(\theta,s) - \left (\partial_\omega + \mu s \frac{\partial}{\partial s} \right ) D_{\theta}W(\theta,s) = O(E), \\
& DX(W(\theta,s)) \partial_s W(\theta,s) - \left (\partial_\omega + \mu s
\frac{\partial}{\partial s} \right ) \partial_{s} W(\theta,s) =
\mu \partial_s W(\theta,s) + O(E)
\end{split}
\end{equation*}

The vectors  $(J\circ W)^{-1} D_\theta WN$ and $(J\circ W )^{-1}
\partial_s W \tilde N$, where $N$ and $\tilde N$ are
normalization matrices (see \eqref{red-whiskers-elem}), are the
symplectic conjugate vectors of $D_\theta W$ and $\partial_s
W$, respectively.

By the Hamiltonian character of the vector field, we have that
\begin{equation*}
\begin{split}
&(DX\circ W(\theta,s)) ((J\circ W)^{-1} D_\theta W  N)(\theta,s) -
\left (
\partial_\omega + \mu s \frac{\partial}{\partial s} \right )
((J\circ W)^{-1} D_\theta W  N ) (\theta,s)= \\
&\qquad\qquad   S(\theta,s) D_\theta W(\theta,s) + O(E)\\
&(DX\circ W(\theta,s)) ((J\circ W)^{-1} \partial_s W \tilde N)
(\theta,s) - \left ( \partial_\omega + \mu s
\frac{\partial}{\partial s} \right ) ((J\circ W)^{-1} \partial_s W
\tilde N)
(\theta,s)=\\
&\qquad\qquad -\mu ((J\circ W)^{-1} \partial_s W \tilde N)
(\theta,s) + B(\theta,s) \partial_s W(\theta,s) + O(E)
\end{split}
\end{equation*}
where $S(\theta,s)$ and $B(\theta,s)$ are matrices which can be
computed.

We summarize these properties in the following proposition:

\begin{proposition}\label{reducibility-whiskers-flow}
Using the same notations \eqref{red-whiskers-elem} as in Proposition
\ref{reducibility-whiskers} and considering the matrix $M(\theta,s)$
introduced in \eqref{M-whiskers}, we have
\begin{eqnarray*}
DX\circ W (\theta,s) M(\theta,s) - \left ( \partial_\omega + \mu
s \frac{\partial}{\partial s} \right ) M(\theta,s) = M(\theta, s) \R
(\theta,s) + O(E)
\end{eqnarray*}
\vskip6pt \noindent where
\begin{equation}\label{matrixtw-flows}
\R(\theta,s)=\left(
\begin{picture}(150,50)
\put(10,20){$\begin{matrix} 0 &S(\theta,s)\\ \noalign{\vskip6pt} 0&
0\end{matrix}$} \put(25,-30){$\bigcirc$} \put(110,20){$\bigcirc$}
\put(85,-30){$\begin{matrix} \mu&B(\theta,s)\\
\noalign{\vskip6pt} 0 &-\mu
\end{matrix}$} \put(0,0){\line(3,0){150}}
\put(75,40){\line(0,-3){90}}
\end{picture}
\right)
\end{equation}
with
\begin{equation*}
\begin{split}
& S(\theta,s)= P^\perp(\theta,s) [ \partial_{\omega} \gamma (\theta,s)
-DX(W(\theta,s)) \gamma (\theta,s)]
\\
& B(\theta,s)= Q^\perp(\theta,s) [ \partial_{\omega} \eta (\theta,s)
-DX(W(\theta,s)) \eta (\theta,s)]\\
\end{split}
\end{equation*}
or, equivalently,
\begin{equation*}
\begin{split}
& S(\theta) = P^\perp(\theta,s) ( \Id - P (\theta,s) \alpha(\theta,s)^\perp
) (DX(W(\theta,s)) + DX(W(\theta,s))^\perp )
P(\theta,s) \\
& B(\theta) = Q^\perp(\theta,s) ( \Id - Q (\theta,s) \beta(\theta,s)^\perp )
(DX(W(\theta,s)) + DX(W(\theta,s))^\perp )
Q(\theta,s) \\
\end{split}
\end{equation*}
and $E$ is the error in \eqref{error-whiskers}.
\end{proposition}

Now if we change the unknowns $\Delta= MV$ in
\eqref{eq:Newton-whiskers-flow} and multiply by $M^{-1}
(\theta,s)$ the LHS, by
Proposition~\ref{reducibility-whiskers-flow}, we are left with the
following system of equations
\begin{equation}\label{Newton-whiskers-mod-flow}
\R (\theta,s) V(\theta,s) - \left (\partial_\omega + \mu s
\frac{\partial}{\partial s} \right ) V(\theta,s) - C(\theta,s)
\delta_u = - \widetilde E(\theta,s) + s\delta H(\theta,s)
\end{equation}
where $\R (\theta,s)$ is given in \eqref{matrixtw-flows} and
\begin{equation*}
\begin{split}
&C (\theta,s) = M^{-1} (\theta,s) (J^{-1}DX) \circ K_0 (\theta),\\
&\widetilde E(\theta,s) = M^{-1}(\theta, s)E(\theta,s),\\
&H(\theta,s) = M^{-1} (\theta, s) \partial_s W (\theta ,s).
\end{split}
\end{equation*}

We expand the terms in \eqref{Newton-whiskers-mod-flow} as a
power series up to some order $L$ and match coefficients of the same
order on both sides of the equation. We use subindices to denote
coordinates and superindices to denote the order. Hence, for order
$s^0$ we have
\begin{align}
& - \partial_{\omega} V_1^0 (\theta) + S^{0}(\theta) V_2^0
(\theta) - C_1^0 (\theta) \delta_\mu = - \widetilde E_1^0 (\theta),
\label{order0-1-flow}\\
&-\partial_{\omega} V_2^0 (\theta)  - C_2^0 (\theta)
\delta_\mu
= -\widetilde E_2^0 (\theta), \label{order0-2-flow}\\
&\mu V_3^0 (\theta) - \partial_{\omega} V_3^0 (\theta)
 +B^0(\theta) V_4^0(\theta) - C_3^0 (\theta)
\delta_\mu = - \widetilde E_3^0 (\theta),
\label{order0-3-flow}\\
& - \mu V_4^0 (\theta)  - \partial_{\omega} V_4^0 (\theta) -
C_4^0 (\theta) \delta_\mu = - \widetilde E_4^0 (\theta).
\label{order0-4-flow}
\end{align}

Notice that \eqref{order0-1-flow} and \eqref{order0-2-flow} can be
solved using the Algorithm \ref{Fast-algorithm-flow}. Hence, we
determine $V_1^0$, $V_2^0$ and $\delta_\mu$. Once we know
$\delta_\mu$, we can solve uniquely for $V_3^0$ and $V_4^0$ in
equations \eqref{order0-3-flow} and \eqref{order0-4-flow}. These
equations do not have any small divisors nor obstructions.

For order $s^1$ we have
\begin{align}
&  - \partial_\omega  V_1^1 (\theta)-
\mu V_1^{1} (\theta) + S^0(\theta) V_2^1(\theta) + S^1 (\theta) V_2^0 (\theta)\label{order1-1-flow}\\
&\qquad \qquad = -\widetilde E_1^1 (\theta) +\delta H_1^0(\theta)
+\delta_\mu C_1^1(\theta),
\notag\\
& - \partial_\omega  V_2^1 (\theta) - \mu V_2^1 (\theta) = -
\widetilde E_2^1 (\theta) + \delta H_2^0 (\theta)
+ \delta_\mu C_2^1 (\theta), \label{order1-2-flow}\\
& - \partial_\omega  V_3^1 (\theta) + B^0 (\theta)
V_4^1(\theta) + B^1 (\theta) V_4^0 (\theta)\label{order1-3-flow}\\
&\qquad\qquad =-\widetilde E_3^1 (\theta) + \delta H_3^0 (\theta) +
\delta_\mu C_3^1(\theta),
\notag\\
&  - \partial_\omega  V_4^1 (\theta) - 2 \mu V_4^1 (\theta) =
-\widetilde E_4^1 (\theta) +\delta H_4^0 (\theta) +\delta_\mu C_4^1
(\theta). \label{order1-4-flow}
\end{align}

Notice that once we choose $\delta$, equations \eqref{order1-1-flow}
and \eqref{order1-2-flow} are uniquely solvable for $V_1^1$ and
$V_2^1$. Recall that $\delta_\mu$ is known, since it has been
computed in the case of order~0 equations.

Similarly, equation~\eqref{order1-4-flow} can be solved without
small divisors nor obstructions. However,
equation~\eqref{order1-3-flow} does have obstructions and small
divisors. In order to overcome this problem, we denote by $F$ and
$G$ the solutions of

\begin{gather*}
- \partial_\omega  F (\theta) - 2 \mu F(\theta) = H_4^0 (\theta),\\
- \partial_\omega G (\theta) - 2 \mu G(\theta) = D_4^1 (\theta)
\end{gather*}
where
$$D_4^1 (\theta) = - \widetilde E_4^1 (\theta) +\delta_\mu C_4^1(\theta),$$
then
$$V_4^1 (\theta) = \delta F(\theta) + G(\theta).$$
Taking   averages  of the equation for $V_3^1$ we get
$$\overline{D_3^1} + \delta \overline{H_3^0} - \overline{B^0F} \delta
- \overline{B^0 G} - \overline{B^1 V_4^0} = 0.$$ So we can solve for
$\delta$ provided that $\overline{H_3^0} - \overline{B^0F} \ne0$.

Now the other orders do not have any obstructions,
\begin{equation}
\begin{split}
& - \partial_\omega  V_1^n (\theta) - n \mu V_1^n (\theta ) +
\sum_{k=0}^n S^{n-k} (\theta) V_2^k (\theta) = - \widetilde E_1^n
(\theta) + \delta H_1^{n-1} (\theta) +
\delta_\mu C_1^n (\theta),\\
\noalign{\vskip6pt}
&- \partial_\omega V_2^n (\theta) - n \mu V_2^n (\theta ) = -
\widetilde E_2^n (\theta) + \delta H_2^{n-1}(\theta)
+ \delta_\mu C_2^n (\theta), \\
\noalign{\vskip6pt}
& -  \partial_\omega  V_3^n (\theta) - (n-1) \mu V_3^n (\theta )
+ \sum_{k=0}^n B^{n-k}(\theta) V_4^k (\theta) = - \widetilde E_3^n
(\theta) + \delta H_3^{n-1}(\theta)
+ \delta_\mu C_3^n (\theta),\\
\noalign{\vskip6pt}
& - \partial_\omega V_4^n (\theta) - (n+1) \mu V_4^n (\theta ) =
- \widetilde E_4^n (\theta) + \delta H_4^{n-1} (\theta) + \delta_\mu
C_4^n (\theta).
\end{split}
\end{equation}
for $n\ge 2$ and they can be solved uniquely
for $V_1^n$, $V_2^n$, $V_3^n$ and $V_4^n$, for $n=2,\ldots , L$, where
$L$ is the degree for the Taylor expansion. Hence, we have
obtained $\delta,\delta_\mu \in \real$ and
$$V (\theta,s) = \sum_{n=0}^L V^n (\theta) s^n$$
and the improved solution is
\begin{equation*}
\begin{split}
& W\gets W+MV\\
& \lambda \gets \lambda +\delta\\
& \mu \gets \mu +\delta_\mu
\end{split}
\end{equation*}

The algorithm for the whiskers and the invariant torus, analoguous to
Algorithm \ref{alg:wt-maps}, is

\begin{algorithm}[Computation of whiskers and tori for flows]
Consider given $X$, $\omega$, $K_0$ and a fixed order $L$. Given an
approximate solution $(W,\mu,\mu)$, perform the following
calculations
\begin{itemize}
\item[1.] Compute $ E(\theta,s)=X(W(\theta,s)) - ( \partial_\omega + \mu s
\partial_s  ) W(\theta,s) - (J^{-1}
DX) \circ K_0(\theta) \mu$
\item[2.] Compute
\begin{itemize}
\item[(2.1)] $\alpha (\theta,s) = D_\theta W(\theta,s)$
\item[(2.2)] $\beta (\theta,s) = \partial_s W (\theta,s)$
\item[(2.3)] $N(\theta,s) = [\alpha (\theta,s)^\perp \alpha (\theta,s)]^{-1}$
\item[(2.4)] $\tilde N(\theta,s) = [\beta (\theta,s)^\perp \beta (\theta,s)]^{-1}$
\item[(2.5)] $P(\theta,s) = \alpha (\theta,s) N (\theta,s)$
\item[(2.6)] $Q(\theta,s) = \beta (\theta,s) \tilde N (\theta,s)$
\item[(2.7)] $\gamma (\theta,s) = (J\circ W (\theta,s))^{-1} P(\theta,s)$
\item[(2.8)] $\eta (\theta,s) = (J\circ W (\theta,s))^{-1} Q(\theta,s)$
\item[(2.9)] $M(\theta,s) = [\alpha (\theta,s) \mid \gamma (\theta,s) \mid
\beta (\theta,s) \mid \eta (\theta,s)]$
\item[(2.10)] $[M(\theta,s)]^{-1}$
\end{itemize}
\item[3.]  Compute
\begin{itemize}
\item[(3.1)] $C(\theta,s) = M^{-1} (\theta,s) (J^{-1}DX) \circ K_0 (\theta)$
\item[(3.2)] $\widetilde E(\theta,s) = M^{-1}(\theta, s)E(\theta,s)$
\item[(3.3)] $H(\theta,s) = M^{-1} (\theta , s) \beta
(\theta ,s)$
\end{itemize}
\item[4.] Compute
\begin{itemize}
\item[(4.1)] $S(\theta,s) = P^\perp(\theta,s) ( \Id - P (\theta,s)
\alpha(\theta,s)^\perp ) (DX(W(\theta,s)) + DX(W(\theta,s))^\perp )
P(\theta,s)$
\item[(4.2)] $B(\theta,s) = Q^\perp(\theta,s) ( \Id - Q (\theta,s)
\beta(\theta,s)^\perp ) (DX(W(\theta,s)) + DX(W(\theta,s))^\perp )
Q(\theta,s)$
\end{itemize}
\item[5.]
\begin{itemize}
\item[(5.1)] Solve for $\delta_\mu$ satisfying
$$\int_{\torus^\ell} \widetilde E_2^0 - \bigg[ \int_{\torus^\ell} C_2^0\bigg] \delta_\mu =0$$
\item[(5.2)] Solve for $V_2^0$ satisfying
$$ - \partial_{\omega} V_2^0  = -\widetilde E_2^0 + C_2^0 \delta_\mu$$
Set $V_2^0$ such that the average is $0$.
\end{itemize}
\item[6.]
\begin{itemize}
\item[(6.1)] Compute $S^0 (\theta) V_2^0 (\theta)$
\item[(6.2)] Solve for $\bar V_2^0$ satisfying
$$\int_{\torus^\ell} \widetilde E_1^0 - \int_{\torus^\ell} C_1^0 (\theta)\delta_\mu
+ \int_{\torus^\ell} S^0 V_2^0 + \bigg[ \int_{\torus^\ell} S^0\bigg] \bar V_2^0=0$$
\item[(6.3)] Set $V_2^0 = V_2^0 + \bar V_2^0$
\item[(6.4)] Solve for $V_1^0$ satisfying
$$- \partial_{\omega} V_1^0  = - \widetilde E_1^0 - S^0 V_2^0
+ C_1^0 \delta_\mu$$
\item[(6.5)] Normalize so that $\int_{\torus^\ell} V_1^0 =0$
\end{itemize}
\item[7.] Solve for $V_4^0$ satisfying
$$ - \mu V_4^0 - \partial_{\omega} V_4^0  =
- \widetilde E_4^0 + C_4^0 \delta_\mu$$
\item[8.] Solve for $V_3^0$ satisfying
$$\mu V_3^0 - \partial_{\omega} V_3^0  = - \widetilde E_3^0
+ C_3^0 \delta_\mu - B^0 V_4^0$$
\item[9.]
\begin{itemize}
\item[(9.1)] Solve for $F$ satisfying
$$ - \partial_{\omega} F - 2 \mu F = H_4^0$$
\item[(9.2)] Solve for $G$ satisfying
$$ - \partial_{\omega} G - 2 \mu G = -\widetilde E_4^1
+\delta_\mu C_4^1$$
\item[(9.3)] Solve for $\delta$ satisfying
$$\left( - \overline{\widetilde E_3^1} + \delta_\mu \overline{C_3^1}
- \overline{B^0 G} - \overline{B^1 V_4^0} \right) +\delta
(\overline{H_3^0} - \overline{B^0 F}) =0$$
\item[(9.4)] Set $V_4^1 = \delta F+G$
\end{itemize}
\item[10.\ \ ]
\begin{itemize}
\item[(10.1)] Solve for $V_3^1$ satisfying
$$ - \partial_{\omega} V_3^1 = -\widetilde E_3^1
+ \delta H_3^0 + \delta_\mu C_3^1 - B^0 V_4^1 - B^1 V_4^0$$
\item[(10.2)] Normalize so that $\int_{\torus^ell} V_3^1=0$
\item[(10.3)] Solve for $V_2^1$ satisfying
$$ - \partial_{\omega} V_2^1 - \mu V_2^1  = - \widetilde E_2^1 + \delta H_2^0
+\delta_\mu C_2^1$$
\item[(10.4)] Solve for $V_1^1$ satisfying
$$- \partial_{\omega} V_1^1 - \mu V_1^1  = -\widetilde E_1^1
+\delta H_1^0 + \delta_\mu C_1^1 - S^0 V_2^1 - S^1 V_2^0$$
\end{itemize}
\item[11.] For $n=2\ldots L$ do
\begin{itemize}
\item[(11.1)] Solve for $V_2^n$ satisfying
$$ - \partial_{\omega} V_2^n - n \mu V_2^n  = -\widetilde E_2^n (\theta)
+ \delta H_2^{n-1} + \delta_\mu C_2^n$$
\item[(11.2)] Compute
$$\tilde S^n = \sum_{k=0}^n S^{n-k} V_2^k$$
\item[(11.3)] Solve for $V_1^n$ satisfying
$$- \partial_{\omega}V_1^n - n \mu V_1^n  = - \widetilde E_1^n
+ \delta H_1^{n-1} + \delta_\mu C_1^n - \tilde S^n$$
\item[(11.4)] Solve for $V_4^n$ satisfying
$$ - \partial_{\omega} V_4^n - (n+1) \mu V_4^n
= - \widetilde E_4^n + \delta H_4^{n-1} + \delta_\mu C_4^n$$
\item[(11.5)] Compute
$$\tilde B^n = \sum_{k=0}^n B^{n-k} V_4^k$$
\item[(11.6)] Solve for $V_3^n$ satisfying
$$ - \partial_{\omega} V_3^n - (n-1) \mu  V_3^n
= - \widetilde E_3^n + \delta H_3^{n-1} + \delta_\mu C_3^n - \tilde
B^n$$
\end{itemize}
\item[12.] Compute
$$V(\theta) = \sum_{n=0}^L V^n (\theta) s^n$$
\item[13.] Set \quad $W\gets W+MV$
\item[{}] \qquad\quad $\lambda \gets \lambda +\delta$
\item[{}] \qquad\quad $\mu \gets \mu +\delta_\mu$
\end{itemize}
\end{algorithm}


\section{A Newton method to compute the whiskers for
flows}\label{ap:d}

We consider the invariance equation~\eqref{invariance-whiskers} for
flows and we assume that we have an initial approximation $W$ for
the whiskers, expressed as a power series
$$W(\theta,s) = \sum_{n=0}^\infty W^n (\theta) s^n$$
and such that
$$W^0 (\theta) = K(\theta)\ \text{ and }\ W^1 (\theta) = V^s(\theta)$$
(the case unstable is analogous).

Then, it is clear that the error $E$ for the initial approximation
$W$ is such that
$$E(\theta,s) = \sum_{n\ge2} E^n (\theta) s^n$$
because the approximation is exact for the terms of order~0 and 1.

Using the notation introduced in \eqref{eq:trunc}, the linearized
equation for the Newton method is
$$[DX \circ W(\theta,s)] \Delta^{[\geq 2 ]} (\theta,s) - ( \partial_\omega + \mu s
\partial s ) \Delta^{[\ge2]}
(\theta,s) = - E^{[\ge2]} (\theta,s).$$

Proceeding as in the previous section we can perform the change of
coordinates given by the matrix $M(\theta,s)$ in \eqref{M-whiskers}
and reduce the problem to solving for $V(\theta,s)$ the following
equation, which is diagonal in Fourier-Taylor series,
$$\R(\theta,s) V^{[\ge 2]} (\theta,s) - ( \partial_\omega + \mu s
\partial s ) V^{[\ge 2]}(\theta, s)
= - \widetilde E^{[\ge 2]} (\theta,s),$$ with $\R(\theta,s)$ given
in \eqref{matrixtw-flows} and $\widetilde E(\theta,s)=M(\theta ,
 s)^{-1} E(\theta,s)$.

Notice that in this case, we do not have to solve the system of
equations for order~0 and 1 and we can go straight to order $n\ge2$.
We use subindices to denote coordinates and superindices to denote
the order. Hence, for order $n \ge 2$, we need to solve the system of
equations
\begin{equation}\label{ordern-simple-flow}
\begin{split}
& - \partial_\omega V_1^n (\theta) - n \mu V_1^n (\theta) +
\sum_{k=2}^n S^{n-k} (\theta) V_2^k (\theta)
= - \widetilde E_1^n (\theta),\\
\noalign{\vskip6pt} &- \partial_{\omega} V_2^n (\theta) - n \mu
V_2^n (\theta)
= - \widetilde E_2^n (\theta),\\
\noalign{\vskip6pt} & - \partial_{\omega} V_3^n (\theta) - (n-1)
\mu V_3^n (\theta
) + \sum_{k=2}^n B^{n-k} (\theta) V_4^k (\theta) = - \widetilde E_3^n,\\
\noalign{\vskip6pt} &  - \partial_{\omega} V_4^n (\theta) - (n+1)
\mu V_1^n (\theta) = - \widetilde E_4^n.
\end{split}
\end{equation}

Notice that now the solution of \eqref{ordern-simple-flow} for
$n=2,3$ provides an exact solution of the invariance equation up to
order~4. That is, if we set
$$V^{[<4]} (\theta,s) = V^2 (\theta,s) + V^3 (\theta,s)$$
where $V^2$ and $V^3$ are obtained by solving
equations~\eqref{ordern-simple-flow}, then the improved solution
$\bar W$ given by
$$ \bar W(\theta,s) = W(\theta,s) + M(\theta,s)V^{[<4]} (\theta,s),$$
where $M(\theta,s)$ was introduced in \eqref{M-whiskers}, satisfies
that it approximates the solution of the invariance equations with
an error $\bar E$ such that
$$ \bar E(\theta,s) = \bar E^{[\ge 4]} (\theta,s).$$

This process can be iterated and at each step we solve the
invariance equation exactly up to an order which is the double of
the one we had for the initial approximation. Thus, if we assume
that the initial guess $W$ is such that the error in
\eqref{error-whiskers-flow} satisfies that
$$E = E^{[\ge L]},$$
then the modified linearized equation for the Newton method is such
that
$$\R(\theta,s) V^{[\ge L]} (\theta,s) - (\partial_{\omega} + \mu s \partial_{s}) V^{[\ge L]} (\theta,s)
= - \widetilde E^{[\ge L]} (\theta,s),$$ with $\R(\theta,s)$ given
in \eqref{matrixtw-flows}. If we solve the system of equations
\eqref{ordern-simple-flow} for $n= L\ldots (2L-1)$ then the improved
$\bar W$ is
$$\bar W(\theta,s) = W(\theta,s) + M(\theta,s) V^{[< 2L]} (\theta,s),$$
with $M(\theta,s)$ as in \eqref{M-whiskers}, and the new error $\bar
E$ satisfies $\bar E(\theta,s)= \bar E^{[\ge 2L]} (\theta,s)$.

The algorithm in this case  is

\begin{algorithm}[Computation of whiskers for vector-fields]
Given $X$, $\omega$ as well as $K, V^s, \mu$ and an approximate
solution $W$ such that
\[ X\circ W(\theta,s) -
(\partial_{\omega} + \mu s \partial_{s}) W(\theta, s) = E^{[\ge
L]}(\theta,s) \] with $L \geq 2$ and $W(\theta,0)=K(\theta)$ and
$\partial_{s}W(\theta,0)=V^s(\theta)$, perform the following
calculations:
\begin{itemize}
\item[1.]
Compute $E^{[\ge L]}(\theta,s) = X\circ W(\theta,s) -
(\partial_{\omega} + \mu s \partial_{s}) W(\theta, s)$
\item[2.] Compute
\begin{itemize}
\item[(2.1)] $\alpha (\theta,s) = D_\theta W(\theta,s)$
\item[(2.2)] $\beta (\theta,s) = \partial_s W (\theta,s)$
\item[(2.3)] $N(\theta,s) = [\alpha (\theta,s)^\perp \alpha (\theta,s)]^{-1}$
\item[(2.4)] $\tilde N(\theta,s) = [\beta (\theta,s)^\perp \beta (\theta,s)]^{-1}$
\item[(2.5)] $P(\theta,s) = \alpha (\theta,s) N (\theta,s)$
\item[(2.6)] $Q(\theta,s) = \beta (\theta,s) \tilde N (\theta,s)$
\item[(2.7)] $\gamma (\theta,s) = (J\circ W (\theta,s))^{-1} P(\theta,s)$
\item[(2.8)] $\eta (\theta,s) = (J\circ W (\theta,s))^{-1} Q(\theta,s)$
\item[(2.9)] $M(\theta,s) = [\alpha (\theta,s) \mid \gamma (\theta,s) \mid
\beta (\theta,s) \mid \eta (\theta,s)]$
\item[(2.10)] $[M(\theta,s)]^{-1}$
\end{itemize}
\item[3.] Compute
$$\widetilde E^{[\ge L]} (\theta,s) = M^{-1} (\theta,s)
E^{[\ge L]} (\theta,s)$$
\item[4.] Compute
\begin{itemize}
\item[(4.1)] $S(\theta,s) = P^\perp(\theta,s) ( \Id - P (\theta,s)
\alpha(\theta,s)^\perp ) (DX(W(\theta,s)) + DX(W(\theta,s))^\perp )
P(\theta,s)$
\item[(4.2)] $B(\theta,s) = Q^\perp(\theta,s) ( \Id - Q (\theta,s)
\beta(\theta,s)^\perp ) (DX(W(\theta,s)) + DX(W(\theta,s))^\perp )
Q(\theta,s)$
\end{itemize}
\item[5.] For $n=L\ldots 2L-1$ do
\begin{itemize}
\item[(5.1)] Solve for $V_2^n$ satisfying
$$- \partial_{\omega} V_2^n (\theta) - n
\mu V_2^n (\theta) = -\widetilde E_2^n (\theta)$$
\item[(5.2)] Compute
$$\tilde S^n = \sum_{k=L}^n S^{n-k} V_2^k$$
\item[(5.3)] Solve for $V_1^n$ satisfying
$$-\partial_{\omega} V_1^n (\theta) - n
\mu V_1^n (\theta) = -\widetilde E_1^n - \tilde S^n$$
\item[(5.4)] Solve for $V_4^n$ satisfying
$$ - \partial_{\omega} V_4^n (\theta) - (n+1)
\mu V_4^n (\theta) = - \widetilde E_4^n$$
\item[(5.5)] Compute
$$\tilde B^n = \sum_{k=L}^n B^{n-k} V_4^k$$
\item[(5.6)] Solve for $V_3^n$ satisfying
$$ - \partial_\omega V_3^n (\theta) - (n-1)
\mu V_3^n (\theta) = - \widetilde E_3^n - \tilde B^n$$
\end{itemize}
\item[6.] Compute
$$V(\theta,s) = \sum_{n=L}^{2L-1} V^n (\theta) s^n$$
\item[7.] Set $W\gets W+MV$
\end{itemize}
\end{algorithm}


\section*{Acknowledgements}
The work of
G.~H. and   R.~L. has been partially supported by NSF grants. G.~H. has
also been supported by the Spanish Grant MTM2006-00478 and the
Spanish Fellowship AP2003-3411.

We thank R. Calleja \'A Haro, A. Luque, J. M. Mondelo and C. Sim\'o 
for several discussions and for
comments on the paper. The final version was written while we were
visiting CRM during the Research Programme \emph{Stability and
Instability in Mechanical Systems (SIMS08)}, for whose hospitality
we are very grateful. G.~H. and  Y.~S.  would like to thank the hospitality of the department of Mathematics of
University of Texas at Austin, where part of this work was carried out.
\bibliographystyle{alpha}

\bibliography{referencies2,new}

\def\cprime{$'$} \def\cprime{$'$} \def\cprime{$'$} \def\cprime{$'$}
  \def\cprime{$'$}
\begin{thebibliography}{FdlLS09b}

\bibitem[Arn64]{Arnold64}
V.I. Arnold.
\newblock Instability of dynamical systems with several degrees of freedom.
\newblock {\em Sov. Math. Doklady}, 5:581--585, 1964.

\bibitem[CC07]{CellettiC07}
Alessandra Celletti and Luigi Chierchia.
\newblock K{AM} stability and celestial mechanics.
\newblock {\em Mem. Amer. Math. Soc.}, 187(878):viii+134, 2007.

\bibitem[CdlL09]{CallejaL08}
Renato Calleja and Rafael de~la Llave.
\newblock Fast numerical computation of quasi-periodic equilibrium states in
  1{D} statistical mechanics, including twist maps.
\newblock {\em Nonlinearity}, 22(6):1311--1336, 2009.

\bibitem[CFL03a]{CabreFL03a}
Xavier Cabr{\'e}, Ernest Fontich, and Rafael de~la Llave.
\newblock The parameterization method for invariant manifolds. {I}. {M}anifolds
  associated to non-resonant subspaces.
\newblock {\em Indiana Univ. Math. J.}, 52(2):283--328, 2003.

\bibitem[CFL03b]{CabreFL03b}
Xavier Cabr{\'e}, Ernest Fontich, and Rafael de~la Llave.
\newblock The parameterization method for invariant manifolds. {II}.
  {R}egularity with respect to parameters.
\newblock {\em Indiana Univ. Math. J.}, 52(2):329--360, 2003.

\bibitem[CFL05]{CabreFL05}
Xavier Cabr{\'e}, Ernest Fontich, and Rafael de~la Llave.
\newblock The parameterization method for invariant manifolds. {III}.
  {O}verview and applications.
\newblock {\em J. Differential Equations}, 218(2):444--515, 2005.

\bibitem[DH09]{DelshamsH09}
Amadeu Delshams and Gemma Huguet.
\newblock Geography of resonances and {A}rnold diffusion in a priori unstable
  {H}amiltonian systems.
\newblock {\em Nonlinearity}, 22(8):1997--2077, 2009.

\bibitem[DLS06]{DelshamsLS06}
A.~Delshams, R.~de~la Llave, and T.~M. Seara.
\newblock A geometric mechanism for diffusion in {H}amiltonian systems
  overcoming the large gap problem: heuristics and rigorous verification on a
  model.
\newblock {\em Mem. Amer. Math. Soc.}, 179(844):viii+141, 2006.

\bibitem[Dua94]{Duarte94}
Pedro Duarte.
\newblock Plenty of elliptic islands for the standard family of area preserving
  maps.
\newblock {\em Ann. Inst. H. Poincar\'e Anal. Non Lin\'eaire}, 11(4):359--409,
  1994.

\bibitem[FdlLS09a]{FontichLS09b}
Ernest Fontich, Rafael de~la Llave, and Yannick Sire.
\newblock Construction of invariant whiskered tori by a parameterization
  method. {I}. {M}aps and flows in finite dimensions.
\newblock {\em J. Differential Equations}, 246(8):3136--3213, 2009.

\bibitem[FdlLS09b]{FontichLS09a}
Ernest Fontich, Rafael de~la Llave, and Yannick Sire.
\newblock A method for the study of whiskered quasi-periodic and
  almost-periodic solutions in finite and infinite dimensional {H}amiltonian
  systems.
\newblock {\em Electron. Res. Announc. Math. Sci.}, 16:9--22, 2009.

\bibitem[FGB98]{FassoGB98a}
F.~Fass{\`o}, M.~Guzzo, and G.~Benettin.
\newblock Nekhoroshev-stability of elliptic equilibria of hamiltonian systems.
\newblock {\em Comm. Math. Phys.}, 197(2):347--360, 1998.

\bibitem[FKW01]{FayadKW01}
Bassam Fayad, Anatole Katok, and Alistar Windsor.
\newblock Mixed spectrum reparameterizations of linear flows on {${\Bbb T}\sp
  2$}.
\newblock {\em Mosc. Math. J.}, 1(4):521--537, 644, 2001.
\newblock Dedicated to the memory of I. G.\ Petrovskii on the occasion of his
  100th anniversary.

\bibitem[GFB98]{GuzzoFB98b}
M.~Guzzo, F.~Fass{\`o}, and G.~Benettin.
\newblock On the stability of elliptic equilibria.
\newblock {\em Math. Phys. Electron. J.}, 4:Paper 1, 16 pp.\ (electronic),
  1998.

\bibitem[Gra74]{Graff74}
Samuel~M. Graff.
\newblock On the conservation of hyperbolic invariant tori for {H}amiltonian
  systems.
\newblock {\em J. Differential Equations}, 15:1--69, 1974.

\bibitem[Har08]{Haro08}
Alex Haro.
\newblock Automatic differentiation tools in computational dynamical systems.
\newblock Manuscript, 2008.

\bibitem[HdlLS09]{HuguetLS10c}
G.~Huguet, R.~de~la Llave, and Y.~Sire.
\newblock Fast iteration of quasi-periodic cocyles.
\newblock {\em Manuscript}, 2009.

\bibitem[Her83]{Herman83}
Michael-R. Herman.
\newblock {\em Sur les courbes invariantes par les diff\'eomorphismes de
  l'anneau. {V}ol. 1}, volume 103 of {\em Ast\'erisque}.
\newblock Soci\'et\'e Math\'ematique de France, Paris, 1983.
\newblock With an appendix by Albert Fathi, With an English summary.

\bibitem[Her92]{Herman92}
M.-R. Herman.
\newblock On the dynamics of {L}agrangian tori invariant by symplectic
  diffeomorphisms.
\newblock In {\em Progress in Variational Methods in Hamiltonian Systems and
  Elliptic Equations (L'Aquila, 1990)}, pages 92--112. Longman Sci. Tech.,
  Harlow, 1992.

\bibitem[HL00]{HaroL00}
A.~Haro and R.~de~la Llave.
\newblock New mechanisms for lack of equipartion of energy.
\newblock {\em Phys. Rev. Lett.}, 89(7):1859--1862, 2000.

\bibitem[HL06a]{HaroL06c}
{\`A}.~Haro and R.~de~la Llave.
\newblock Manifolds on the verge of a hyperbolicity breakdown.
\newblock {\em Chaos}, 16(1):013120, 8, 2006.

\bibitem[HL06b]{HaroL06b}
{\`A}.~Haro and R.~de~la Llave.
\newblock A parameterization method for the computation of invariant tori and
  their whiskers in quasi-periodic maps: numerical algorithms.
\newblock {\em Discrete Contin. Dyn. Syst. Ser. B}, 6(6):1261--1300
  (electronic), 2006.

\bibitem[HL06c]{HaroL06a}
A.~Haro and R.~de~la Llave.
\newblock A parameterization method for the computation of invariant tori and
  their whiskers in quasi-periodic maps: rigorous results.
\newblock {\em J. Differential Equations}, 228(2):530--579, 2006.

\bibitem[HL07]{HaroL07}
A.~Haro and R.~de~la Llave.
\newblock A parameterization method for the computation of whiskers in quasi
  periodic maps: numerical implementation and examples.
\newblock {\em SIAM Jour. Appl. Dyn. Syst.}, 6(1):142--207, 2007.

\bibitem[JO05]{JorbaO03}
{\`A}ngel Jorba and Estrella Olmedo.
\newblock A parallel method to compute quasi-periodic solutions.
\newblock In {\em E{QUADIFF} 2003}, pages 181--183. World Sci. Publ.,
  Hackensack, NJ, 2005.

\bibitem[JO09]{JorbaO09}
{\`A}ngel Jorba and Estrella Olmedo.
\newblock On the computation of reducible invariant tori on a parallel
  computer.
\newblock {\em SIAM J. Appl. Dyn. Syst.}, 8(4):1382--1404, 2009.

\bibitem[Knu97]{Knuth97}
Donald~E. Knuth.
\newblock {\em The art of computer programming. {V}ol. 2: {S}eminumerical
  algorithms}.
\newblock Addison-Wesley Publishing Co., Reading, Mass.-London-Don Mills, Ont,
  third revised edition, 1997.

\bibitem[LGJV05]{LlaveGJV05}
R.~de~la Llave, A.~Gonz{\'a}lez, {\`A}.~Jorba, and J.~Villanueva.
\newblock K{AM} theory without action-angle variables.
\newblock {\em Nonlinearity}, 18(2):855--895, 2005.

\bibitem[Lla01]{Llave01}
Rafael de~la Llave.
\newblock A tutorial on {K}{A}{M} theory.
\newblock In {\em Smooth ergodic theory and its applications (Seattle, WA,
  1999)}, pages 175--292. Amer. Math. Soc., Providence, RI, 2001.

\bibitem[LW04]{LlaveW04}
R.~de~la Llave and C.~E. Wayne.
\newblock Whiskered and low dimensional tori in nearly integrable {H}amiltonian
  systems.
\newblock {\em Math. Phys. Electron. J.}, 10:Paper 5, 45 pp. (electronic),
  2004.

\bibitem[RC95]{RakovicChu95}
M.~J. Rakovi{\'c} and Shih-I Chu.
\newblock New integrable systems: hydrogen atom in external fields.
\newblock {\em Phys. D}, 81(3):271--279, 1995.

\bibitem[RC97]{RakovicChu97}
M.~J. Rakovi{\'c} and Shih-I Chu.
\newblock Phase-space structure of a new integrable system related to hydrogen
  atoms in external fields.
\newblock {\em J. Phys. A}, 30(2):733--753, 1997.

\bibitem[R{\"u}s75]{Russmann75}
H.~R{\"u}ssmann.
\newblock On optimal estimates for the solutions of linear partial differential
  equations of first order with constant coefficients on the torus.
\newblock In {\em Dynamical Systems, Theory and Applications (Battelle
  Rencontres, Seattle, Wash., 1974)}, pages 598--624. Lecture Notes in Phys.,
  Vol. 38, Berlin, 1975. Springer.

\bibitem[Sim00]{Simo00}
Carles Sim\'o.
\newblock private communication.
\newblock 2000.

\bibitem[Zeh76]{Zehnder76}
E.~Zehnder.
\newblock Generalized implicit function theorems with applications to some
  small divisor problems. {I}{I}.
\newblock {\em Comm. Pure Appl. Math.}, 29(1):49--111, 1976.

\end{thebibliography}

\end{document}